\documentclass[12pt]{article}
\usepackage{amsmath,amsfonts,amssymb,rotating}

\def\qed{\nopagebreak\hfill{\rule{4pt}{7pt}}}
\def\proof{\noindent {\it{Proof.} \hskip 2pt}}
\def\ps{\mathrm{ps}}

\newtheorem{theo}{Theorem}[section]

\newtheorem{lemm}[theo]{Lemma}

\newtheorem{coro}[theo]{Corollary}
\newtheorem{conj}[theo]{Conjecture}

\newdimen\Squaresize \Squaresize=11pt
\newdimen\Thickness \Thickness=0.7pt
\def\Square#1{\hbox{\vrule width \Thickness
   \vbox to \Squaresize{\hrule height \Thickness\vss
    \hbox to \Squaresize{\hss#1\hss}
   \vss\hrule height\Thickness}
\unskip\vrule width \Thickness} \kern-\Thickness}

\def\Vsquare#1{\vbox{\Square{$#1$}}\kern-\Thickness}

\def\young#1{
\vbox{\smallskip\offinterlineskip \halign{&\Vsquare{##}\cr #1}}}

\def\moins{\raise 1pt\hbox{{$\scriptstyle -$}}}

\begin{document}

\parskip 8pt

\begin{center}
{\large \bf Schur Positivity and the  $q$-Log-convexity of

 the
Narayana Polynomials}
\end{center}

\begin{center}
William Y. C. Chen$^{1}$, Larry X.W. Wang$^{2}$ and Arthur
L. B. Yang$^{3}$\\[6pt]
Center for Combinatorics, LPMC-TJKLC\\
Nankai University, Tianjin 300071, P. R. China

 Email: $^{1}${\tt
chen@nankai.edu.cn}, $^{2}${\tt wxw@cfc.nankai.edu.cn}, $^{3}${\tt
yang@nankai.edu.cn}
\end{center}

\vspace{0.3cm} \noindent{\bf Abstract.} Using Schur positivity and
the principal specialization of Schur functions, we provide a proof
of a recent conjecture of Liu and Wang on the $q$-log-convexity of
the Narayana polynomials, and a proof of the second conjecture that
the
 Narayana transformation preserves the
log-convexity. Based on a formula of Br\"and$\mathrm{\acute{e}}$n
which expresses the $q$-Narayana numbers as the specializations of
Schur functions, we derive several symmetric function identities
using the Littlewood-Richardson rule for the product of Schur
functions, and obtain the strong $q$-log-convexity of the Narayana
polynomials and the strong $q$-log-concavity of the $q$-Narayana
numbers.

\noindent {\bf Keywords:} $q$-log-concavity, $q$-log-convexity,
$q$-Narayana number, Narayana polynomial, lattice permutation, Schur
positivity, Littlewood-Richardson rule.

\noindent {\bf AMS Classification:} 05E05, 05E10

\noindent {\bf Suggested Running Title:} Schur positivity

\section{Introduction}

The main objective of this paper is to provide proofs of two recent
conjectures of Liu and Wang \cite{liuwan2006} on the
$q$-log-convexity of the Narayana polynomials by using Schur
positivity derived from the Littlewood-Richardson rule. Moreover, we
prove that the Narayana polynomials are strongly $q$-log-convex. We
also study the $q$-log-concavity of the $q$-Narayana numbers, and
prove that for fixed $n$ or $k$ the $q$-Narayana numbers $N_q(n,k)$
 are strongly $q$-log-concave.

Unimodal and log-concave sequences and polynomials often arise in
 combinatorics, algebra and geometry, see, for example, Brenti \cite{brenti1989, brenti1994}, Stanley \cite{stanley1989},
and Stembridge \cite{stembr2007}. A sequence $(a_n)_{n\geq 0}$ of
real numbers is said to be unimodal if there exists an integer
$m\geq 0$ such that
$$a_0\leq a_1\leq\cdots\leq a_m\geq a_{m+1}\geq
a_{m+2} \geq \cdots,$$ and is said to be log-concave if
$$a_m^2\geq a_{m+1}a_{m-1}$$  holds for all $m\geq 1$.

It has been noticed that sometimes the reciprocals of a
combinatorial sequence form a log-concave sequence. For example, the
sequence
$$\binom{n}{0}^{-1},\binom{n}{1}^{-1},\ldots,\binom{n}{n}^{-1}$$
satisfies this condition for a given positive integer $n$. Such
sequences are called log-convex, see \cite{liuwan2006}.

For polynomials, Stanley introduced the notion of $q$-log-concavity,
which has been studied  by Butler \cite{butler1990}, Krattenthaler
\cite{kratte1989}, Leroux \cite{leroux1990}, and Sagan
\cite{sagan1992}. A sequence of polynomials $(f_n(q))_{n\geq 0}$
over the field of real numbers is called $q$-log-concave if the
difference
$$f_{m}(q)^2-f_{m+1}(q)f_{m-1}(q)$$
has nonnegative coefficients as a polynomial of $q$ for all $m\geq
1$. Sagan \cite{sagan1992t} also introduced the notion of strong
$q$-log-concavity. We say that a sequence of polynomials
$(f_n(q))_{n\geq 0}$ is strongly $q$-log-concave if
$$f_{m}(q)f_{n}(q)-f_{m+1}(q)f_{n-1}(q)$$
has nonnegative coefficients for any $m\geq n\geq 1$.

Based on the $q$-log-concavity, it is natural to define the
$q$-log-convexity. We say that the polynomial sequence
$(f_n(q))_{n\geq 0}$ is $q$-log-convex if the difference
$$f_{m+1}(q)f_{m-1}(q)-f_{m}(q)^2$$
has nonnegative coefficients as a polynomial of $q$ for all $m\geq
1$. The notion of strong $q$-log-convexity is a natural counterpart
of that of strong $q$-log-concavity. We say that a sequence of
polynomials $(f_n(q))_{n\geq 0}$ is strongly $q$-log-convex if
$$f_{m+1}(q)f_{n-1}(q)-f_{m}(q)f_{n}(q)$$
has nonnegative coefficients for any $m\geq n\geq 1$.

As realized by Sagan \cite{sagan1992t},  the strong
$q$-log-concavity is not equivalent to the $q$-log-concavity,
although for a sequence of positive  numbers  the strong
log-concavity is equivalent to the log-concavity. Analogously, the
strong $q$-log-convexity is not equivalent to the $q$-log-convexity.
For example, the sequence
$$2q+q^2+3q^3,q+2q^2+2q^3,q+2q^2+2q^3,2q+q^2+3q^3$$ is
$q$-log-convex, but not strongly $q$-log-convex.

Just recently, Liu and Wang \cite{liuwan2006} have shown that  some
well-known polynomials such as the Bell polynomials and the Eulerian
polynomials are $q$-log-convex, and proposed some conjectures on the
Narayana polynomials  based on numerical evidence.

To describe the conjectures of Liu and Wang \cite{liuwan2006}, we
begin with the classical Catalan numbers, as given by
 \[ C_n=\frac{1}{n+1}\binom{2n}{n},\]
  which count
the number of Dyck paths from $(0,0)$ to $(2n,0)$ with up steps
$(1,1)$ and down steps $(1,-1)$ but never going below the $x$-axis,
see, Stanley \cite{stanley1999}.  It is known that the Catalan
numbers $C_n$ form a log-convex sequence. Recall that a peak of a
Dyck path is defined as a point where an up step  is immediately
followed by a down  step. Then the Narayana number
\[ N(n,k)=\frac{1}{n}\binom{n}{k}\binom{n}{k+1}\]
 equals the
number of Dyck paths of length $2n$ with exactly $k+1$ peaks, see
 \cite{braden2004, deusch1998, sulank1998,
sulank1999}. The Narayana polynomials  are given by
$$N_n(q)=\sum_{k=0}^n N(n,k)q^k.$$ Liu
and Wang \cite{liuwan2006} have shown that for a given positive real
number $q$ the sequence $(N_n(q))_{n\geq 0}$ is log-convex. Note
that the sequence of the Catalan numbers becomes a  special case for
$q=1$.  The first conjecture of Liu and Wang is as follows.

\begin{conj} \label{mainprob}
The Narayana polynomials $N_n(q)$ form a $q$-log-convex sequence.
\end{conj}

We will prove the above conjecture by studying the Schur positivity
of certain sums of symmetric functions. Our proof heavily relies on
the Littlewood-Richardson rule for the product of Schur functions of
certain shapes with only two columns.  It is the formula of
Br\"and$\mathrm{\acute{e}}$n \cite{braden2004} that enables us to
represent the Narayana polynomials in terms of Schur functions.

To prove the desired Schur positivity, we need to verify several
identities on Schur functions, and we would acknowledge the powerful
role of the Maple packages for symmetric functions, ACE
\cite{veigneau1998} and SF \cite{stembridge1995}.

The second conjecture of Liu and Wang \cite{liuwan2006} is concerned
with the Narayana transformation on sequences of positive real
numbers. The Davenport-P$\mathrm{\acute{o}}$lya theorem
\cite{davpol1949} states that if $(a_n)_{n\geq 0}$ and $(b_n)_{n\geq
0}$ are log-convex then their binomial convolution
$$c_n=\sum_{k=0}^n\binom{n}{k}a_kb_{n-k},\quad n\geq 0$$
is also log-convex. It is known that the binomial convolution also
preserves the log-concavity \cite{wanyeh2007}. However, it is not
generally true that a log-convexity preserving transformation
  also preserves the log-concavity. The componentwise sum is  a simple
  example. On the other hand, Liu and Wang \cite{liuwan2006} have realized that the componentwise sum
  preserves log-convexity.  Moreover, it is  not true that a transformation
which preserves  log-concavity necessarily preserves log-convexity.
The ordinary convolution is such an example \cite{liuwan2006,
wanyeh2007}.

Given combinatorial numbers $(t(n,k))_{0\leq k\leq n}$ such as the
binomial coefficients, one can define a linear operator which
transforms a sequence $(a_n)_{n\geq 0}$ into another sequence
$(b_n)_{n\geq 0}$ as given by
$$b_n=\sum_{k=0}^n t(n,k)a_k, \quad n\geq 0.$$
Liu and Wang \cite{liuwan2006} have  shown that the log-convexity is
preserved by linear transformations associated with  the binomial
coefficients, the Stirling numbers of the first kind and the second
kind. The following conjecture is due to Liu and Wang
\cite{liuwan2006}.

\begin{conj}\label{conj-2}
The Narayana transformation $b_n=\sum_{k=0}^n N(n,k)a_k$ preserves
log-convexity.
\end{conj}

We will give a proof of this conjecture based on the monotone
property of certain quartic polynomials and the $q$-log-convexity of
Narayana polynomials.

In addition, we further prove the strong $q$-log-concavity of the
$q$-Narayana numbers. The $q$-Narayana numbers,  as a natural
$q$-analogue of the Narayana numbers $N(n,k)$, arise from the study
of $q$-Catalan numbers \cite{furhof1985}. The $q$-Narayana number
$N_{q}(n,k)$ is given by
\begin{equation}
N_{q}(n,k)=\frac{1}{[n]}\left[{n\atop k}\right]\left[{n\atop
{k+1}}\right]q^{k^2+k},
\end{equation}
where we use the standard notation
$$[k]:=(1-q^k)/(1-q),\quad [k]!=[1][2]\cdots [k], \quad \left[{n\atop j}\right]:=\frac{[n]!}{[j]![n-j]!}$$
for the $q$-analogues of the integer $k$, the $q$-factorial, and the
$q$-binomial coefficient.

It is known that the $q$-Narayana number $N_{q}(n,k)$ is the natural
refinement of the $q$-Catalan number
$c_n(1)=\frac{1}{[n+1]}\left[{2n\atop {n}}\right]$ defined in
\cite{furhof1985}. Br\"and$\mathrm{\acute{e}}$n \cite{braden2004}
studied several Narayana statistics and bi-statistics on Dyck paths,
and noticed that the $q$-Narayana number $N_{q}(n,k)$ has a Schur
function expression by specializing the variables.

\begin{theo}[{\cite[Theorem 6]{braden2004}}]
For all $n,k\in \mathbb{N}$ we have
\begin{equation}\label{formula-b}
N_{q}(n,k)=s_{(2^k)}(q,q^2,\ldots,q^{n-1}).
\end{equation}
\end{theo}

It was known that the $q$-analogues of many well-known combinatorial
numbers are strongly $q$-log-concave. Bulter \cite{butler1990} and
Krattenthaler \cite{kratte1989} proved the $q$-log-concavity of the
$q$-binomial coefficients, and Leroux \cite{leroux1990} and Sagan
\cite{sagan1992} studied the $q$-log-concavity of the $q$-Stirling
numbers of the first kind and the second kind. It was also known
that the Narayana numbers $N(n,k)$ are log-concave for given $n$ or
$k$. Based on some symmetric function identities, we will show that
$N_{q}(n,k)$ are strongly $q$-log-concave for given $n$ or $k$.

This paper is organized as follows. In Section 2,  we give a brief
review of relevant background on symmetric functions. In Section 3,
we give several symmetric function identities involving  Schur
functions indexed by two-column shapes, and derive the Schur
positivity needed to prove the two conjectures of Liu and Wang.
Section 4 deals with the strong $q$-log-convexity of Narayana
polynomials. The notion of strong $q$-log-convexity is analogous to
that of strong $q$-log-concavity as given by Sagan
\cite{sagan1992t}. In Section 5, we show that  the Narayana
transformation preserves log-convexity. Finally, in Section 6 we
derive  the strong $q$-log-concavity of the $q$-Narayana numbers.

\section{Background on Symmetric Functions}

In this section we  review some relevant background on symmetric
functions and present several recurrence formulas for computing the
principal specializations of Schur functions indexed by certain
two-column shapes, which will be used later in the proofs of the
main theorems. More specifically, the hook-content formula plays an
important role in reducing  the log-convexity preserving property of
the Narayana transformation to the monotone property of certain
polynomials, and the recurrence formulas enable us to reduce the
$q$-log-convexity for Narayana polynomials to the Schur positivity
for certain sums of symmetric functions.

Throughout this paper we will adopt the notation and terminology on
partitions and symmetric functions in Stanley \cite{stanley1999}.
Given a nonnegative integer $n$, a partition $\lambda$ of $n$ is a
weakly decreasing nonnegative integer sequence
$(\lambda_1,\lambda_2,\ldots,\lambda_k)\in\mathbb{N}^k$ such that
$\sum_{i=1}^k\lambda_i=n$. The number of nonzero components
$\lambda_i$ is called the length of $\lambda$, denoted
$\ell(\lambda)$. We also denote the partition $\lambda$ by
$(\ldots,2^{m_2},1^{m_1})$ if $i$ appears $m_i$ times in $\lambda$.
For example, $\lambda=(4,2,2,1,1,1)=(4^1,2^2,1^3)$, where we omit
$i^{m_i}$ if $m_i=0$. Let $\mathrm{Par}(n)$ denote the set of all
partitions of $n$. The Young diagram of $\lambda$ is an array of
squares in the plane justified from the top and left corner with
$\ell(\lambda)$ rows and $\lambda_i$ squares in row $i$. By
transposing the diagram of $\lambda$, we get the conjugate partition
of $\lambda$, denoted $\lambda'$. A square $(i,j)$ in the diagram of
$\lambda$ is the square in row $i$ from the top and column $j$ from
the left. The hook length of $(i,j)$, denoted $h(i,j)$, is given by
$\lambda_i+\lambda_j'-i-j+1$. The content of $(i,j)$, denoted
$c(i,j)$, is given by $j-i$. Given two partitions $\lambda$ and
$\mu$, we say that $\lambda$ contains $\mu$, denoted
$\mu\subseteq\lambda$, if $\lambda_i\geq \mu_i$ holds for each $i$.
When $\mu\subseteq \lambda$, we can define a skew partition
$\lambda/\mu$ as the diagram obtained from the diagram of $\lambda$
by removing the diagram of $\mu$ at the top-left corner.

A semistandard Young tableau of shape $\lambda/\mu$ is an array
$T=(T_{ij})$ of positive integers of shape $\lambda/\mu$ that is
weakly increasing in every row and strictly increasing in every
column. The type of $T$ is defined as the composition
$\alpha=(\alpha_1,\,\alpha_2,\ldots)$, where $\alpha_i$ is the
number of $i$'s in $T$. Let $x$ denote the set of variables $\{
x_1,x_2,\ldots\}$. If $T$ has type $\mathrm{type}(T)=\alpha$, then
we write
$$x^T=x_1^{\alpha_1}x_2^{\alpha_2}\cdots.$$
The skew Schur function $s_{\lambda/\mu}(x)$ of shape $\lambda/\mu$
is defined as the generating function
$$s_{\lambda/\mu}(x)=\sum_T x^T,$$
summed over all semistandard Young tableaux $T$ of shape
$\lambda/\mu$ filled with positive integers. When $\mu$ is the empty
partition $\emptyset$, we call $s_{\lambda}(x)$ the Schur function
of shape $\lambda$. In particular, we set $s_{\emptyset}(x)=1$. It
is well known that the Schur functions $s_{\lambda}$ form a basis
for the ring of symmetric functions.

Let $y=\{y_1, y_2, \ldots\}$ be another set of variables, and  let
$s_{\lambda/\mu}(x, y)$ denote the Schur function in  $x \cup y$.
Note that
\begin{equation}\label{eq-expansion}
s_{\lambda/\mu}(x, y)=\sum_{\nu}s_{\lambda/\nu}(x)s_{\nu/\mu}(y),
\end{equation}
where the sum ranges over all partitions $\nu$ satisfying
$\mu\subseteq \nu\subseteq \lambda$, see \cite{macdonald1995,
stanley1999}.

For a symmetric function $f(x)$, its principle specialization
$\ps_n(f)$ and specialization $\ps_n^1(f)$ of order $n$ are defined
by
$$
\begin{array}{rcl}
\ps_n(f) & = & f(1,q,\ldots,q^{n-1}),\\
\ps_n^1(f) & = & \ps_n(f)|_{q=1}=f(1^n).
\end{array}
$$
For notational convenience, we often omit the variable set $x$ and
simply write $s_{\lambda}$ for the Schur function $s_{\lambda}(x)$
if no confusion arises in the context. The following formula is
called the hook-content formula due to Stanley \cite{stanley1971}.

\begin{lemm}[{\cite[Corollary
7.21.4]{stanley1999}}]\label{lemm-hook}
For any partition $\lambda$ and $n\geq 1$, we have
\begin{equation}
\ps_n(s_{\lambda})=q^{\sum_{k\geq 1}(k-1)\lambda_k}\prod_{(i,j)\in
\lambda}\frac{[n+c(i,j)]}{[h(i,j)]}
\end{equation}
and
\begin{equation}
\ps_n^1(s_{\lambda})=\prod_{(i,j)\in
\lambda}\frac{n+c(i,j)}{h(i,j)}.
\end{equation}
\end{lemm}

On the other hand, in view of \eqref{eq-expansion},  we deduce the
following formulas for the principle specializations of the Schur
functions $s_{\lambda}$ indexed by two-column shapes.

\begin{lemm} \label{lemm-ps} Let $k$ be a positive integer and $n>1$. For any $a<0$ or $b<0$, set $s_{(2^a,1^b)}=0$ by convention.
Then we have
\begin{eqnarray}\label{eq-p}
{\ps_n\left(s_{(2^k)}\right)}=\ps_{n-1}\left(s_{(2^{k})}\right)&+& q^{n-1} \ps_{n-1}\left(s_{(2^{k-1},1)}\right)\nonumber\\
&+& q^{2(n-1)} \ps_{n-1}\left(s_{(2^{k-1})}\right)
\end{eqnarray}
and
\begin{eqnarray}
\ps_n\left(s_{(2^k,1)}\right)=\ps_{n-1}\left(s_{(2^{k},1)}\right)&+& q^{n-1} \ps_{n-1}\left(s_{(2^{k})} + s_{(2^{k-1},1^2)}\right)\nonumber\\
&+& q^{2(n-1)} \ps_{n-1}\left(s_{(2^{k-1},1)}\right).
\end{eqnarray}
Furthermore,
\begin{eqnarray}
\ps_n^1\left(s_{(2^k)}\right) & = &
\ps_{n-1}^1\left(s_{(2^{k})}+s_{(2^{k-1},1)}+s_{(2^{k-1})}\right),\\
\ps_n^1\left(s_{(2^k,1)}\right) & = &
\ps_{n-1}^1\left(s_{(2^{k},1)}+s_{(2^{k})} +
s_{(2^{k-1},1^2)}+s_{(2^{k-1},1)}\right).
\end{eqnarray}
\end{lemm}

\begin{lemm}\label{lemm-general}
For any $m\geq n\geq 1$ and $k\geq 0$, we have
\begin{equation}
\ps_m^1\left(s_{(2^k)}\right)=\sum_{0\leq a\leq b\leq m-n}
\ps_n^1(s_{(2^{k-b},1^{b-a})})\ps_{m-n}^1(s_{(2^a,1^{b-a})}).
\end{equation}
\end{lemm}

The Littlewood-Richardson rule enables us  to expand a product of
Schur functions in terms of Schur functions.  There are several
versions of the Littlewood-Richardson rule; see \cite[Chapter 7,
Appendix A1.3]{stanley1999} and \cite[Part I, Chapter
5]{fulton1997}. These settings have their own advantages when
applied to various problems. For example, Knutson and Tao
\cite{knutao1999} used the honeycomb model to prove the saturation
conjecture. A well-known version is the combinatorial interpretation
of the Littlewood-Richardson coefficients in terms of lattice
permutations, which we will adopt for our purpose.

Recall that a lattice permutation of length $n$ is a sequence
$w_1w_2\cdots w_n$ such that for any $i$ and $j$ in the subsequence
$w_1w_2\cdots w_j$ the number of $i$'s is greater than or equal to
the number of $i+1$'s. Let $T$ be a semistandard Young tableau. The
reverse reading word $T^{\mathrm{rev}}$ is a sequence of entries of
$T$ obtained by first reading each row from right to left and then
concatenating the rows from top to bottom. If the reverse reading
word $T^{\mathrm{rev}}$ is a lattice permutation, we call $T$ a
Littlewood-Richardson tableau. Given two Schur functions $s_{\mu}$
and $s_{\nu}$, Littlewood-Richardson coefficients
$c_{\mu\nu}^\lambda$ can be defined by the following relation
\begin{equation}\label{schur-prod}
s_{\mu}s_{\nu}=\sum_{\lambda}c_{\mu\nu}^{\lambda}s_{\lambda}.
\end{equation}

\begin{theo}[{\cite[Theorem A1.3.3]{stanley1999}}]\label{lrrule}
The Littlewood-Richardson coefficient $c_{\mu\nu}^{\lambda}$ is
equal to the number of Littlewood-Richardson tableaux of shape
$\lambda/\mu$ and type $\nu$.
\end{theo}

Take $\lambda=(9,5,3,3,1),\,\mu=(4,2,1),\,\nu=(7,4,3)$. By using the
Maple package for symmetric functions we find that
$c_{\mu\nu}^{\lambda}=3$. Indeed, there are three
Littlewood-Richardson tableaux of shape $\lambda/\mu$ and type $\nu$
as shown in Figure \ref{fig-1}.

\begin{figure}[h,t]
$$
\young{ * & * & * & * & 1 & 1 & 1 & 1 & 1 \cr
        * & * & 1 & 1 & 2 \cr
        * & 2 & 2 \cr
         2 & 3 & 3 \cr
         3 \cr
         }
         \quad
\young{ * & * & * & * & 1 & 1 & 1 & 1 & 1 \cr
        * & * & 1 & 2 & 2\cr
        * & 2 & 2 \cr
         1 & 3 & 3 \cr
         3 \cr
         }
         \quad
\young{ * & * & * & * & 1 & 1 & 1 & 1 & 1 \cr
        * & * & 1 & 2 & 2 \cr
        * & 1 & 2 \cr
         2 & 3 & 3 \cr
         3 \cr
         }
$$
\caption{Skew Littlewood-Richardson tableaux}\label{fig-1}
\end{figure}

When taking $\nu=(n)$ or $\nu=(1^n)$ in \eqref{schur-prod}, the
Littlewood-Richardson rule has a simpler description, known as
Pieri's rule. We need the notion of horizontal and vertical strips.
A skew partition $\lambda/\mu$ is called a horizontal (or vertical)
strip if there are no two squares in the same column (resp. in the
same row).
\begin{theo}[{\cite[Theorem 7.15.7, Corollary 7.15.9]{stanley1999}}]\label{pierirule}
We have
$$s_{\mu}s_{(n)}=\sum_{\lambda}s_{\lambda}$$
summed over all partitions $\lambda$ such that $\lambda/\mu$ is a
horizontal strip of size $n$, and
$$s_{\mu}s_{(1^n)}=\sum_{\lambda}s_{\lambda}$$
summed over all partitions $\lambda$ such that $\lambda/\mu$ is a
vertical strip of size $n$.
\end{theo}

\section{Schur positivity}

The main goal of this section is to prove the Schur positivity of
certain sums of symmetric functions, which will be needed in the
proof of the $q$-log-convexity of the Narayana polynomials in
Section 4. Given a symmetric function $f$, recall that $f$ is called
$s$-positive (or $s$-negative) if the coefficients $a_\lambda$ in
the expansion $f=\sum_{\lambda}a_{\lambda}s_{\lambda}$ of $f$ in
terms of Schur functions are all nonnegative (resp. nonpositive).

The Schur positivity we establish is deduced from several symmetric
function identities which will be proved by induction based on the
Littlewood-Richardson rule. More specifically, the identities we
consider will involve  the products of Schur functions indexed by
partitions with only two-columns. These Schur functions are of
particular interest for their own sake, see, for example, Rosas
\cite{rosas2001}, and Remmel and Whitehead \cite{remwhi1994}.

It is time to mention that throughout this paper, the
Littlewood-Richardson coefficients are either one or two, and we
will consider those shapes that will occur in the expansion of the
product of Schur functions. It is worth mentioning that the Schur
expansion of the product of two Schur functions would be
multiplicity free when one factor is indexed by a rectangular shape,
see Stembridge \cite{stem2001}.

Let us first introduce certain classes of products of Schur
functions that will be the ingredients to establish the desired
Schur positivity. Given $m\in\mathbb{N}$ and $0\leq i\leq m$, let
\begin{eqnarray*}
D_{m,i}^{(1)} & = &
s_{(2^{i})}s_{(2^{m-i-1})},\\
D_{m,i}^{(2)} & = & s_{(2^{i-1},1^2)}s_{(2^{m-i-1})},\\
D_{m,i}^{(3)} & = & s_{(2^{i-1},1)}s_{(2^{m-i-1},1)},
\end{eqnarray*}
and let
\begin{equation}
D_{m,i}=D_{m,i}^{(1)}+D_{m,i}^{(2)}-D_{m,i}^{(3)},
\end{equation}
where $s_{(2^{i},1)}=s_{(2^{i},1^2)}=0$ for $i<0$ by convention. It
is clear that $D_{m,m}\equiv 0$.

 For two
partitions $\lambda$ and $\mu$, let $\lambda\cup\mu$ be the
partition obtained by taking the union of all parts of $\lambda$ and
$\mu$ and then rearranging them in the weakly decreasing order. For
$k\in \mathbb{N}$ we use $\lambda^k$ to represent the union of $k$
$\lambda$'s, and in particular put $\lambda^k=\emptyset$ if $k=0$.
In this notation, we  introduce an operator $\Delta^{\mu}$ on the
ring of symmetric functions defined by a partition $\mu$. For a
symmetric function $f$, suppose that $f$ has the expansion
$$f=\sum_{\lambda}a_{\lambda}s_{\lambda},$$ and then the action of
$\Delta^{\mu}$ on $f$ is given by
$$\Delta^{\mu}(f)=\sum_{{\lambda}}a_{\lambda}s_{\lambda\cup\mu}.$$
For example, if
$$f=s_{(4,3,2)}+3s_{(2,2,1)}+2s_{(5)},$$
then
$$\Delta^{(3,1)}f=s_{(4,3,3,2,1)}+3s_{(3,2,2,1,1)}+2s_{(5,3,1)}.$$

\begin{lemm} \label{lemm-first}
For any $n\geq k\geq 1$, we have
\begin{eqnarray}
s_{(2^{k})}s_{(2^{n+1})} & = & \Delta^{(2)}(s_{(2^{k})}s_{(2^{n})}), \label{eq-s-1}\\
s_{(2^{k-1},1^2)}s_{(2^{n+1})}  & = & \Delta^{(2)}(s_{(2^{k-1},1^2)}s_{(2^{n})}),\\
s_{(2^{k})}s_{(2^{n+1},1^2)}  & = & \Delta^{(2)}(s_{(2^{k})}s_{(2^{n},1^2)}),\\
s_{(2^{k-1},1)}s_{(2^{n+1},1)}  & = &
\Delta^{(2)}(s_{(2^{k-1},1)}s_{(2^{n},1)}).
\end{eqnarray}
\end{lemm}

\begin{figure}[h,t]
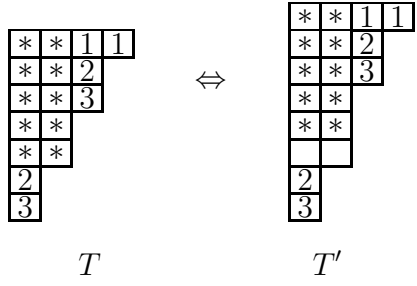

$$
\young{ * & * & 1 & 1\cr
        * & * & 2 \cr
        * & * & 3 \cr
        *  & * \cr
        *  & * \cr
         2\cr
         3\cr
         }
         \quad  \quad \raisebox{50pt}{$\Leftrightarrow$} \quad  \quad
\young{ * & * & 1 & 1\cr
        * & * & 2 \cr
        * & * & 3 \cr
        * & * \cr
        * & * \cr
         &  \cr
         2\cr
         3\cr
         }
$$
$$T \hspace{80pt} T'$$
\caption{Bijection between Littlewood-Richardson
tableaux}\label{fig-2}
\end{figure}

\proof Define $a_\lambda$ by
\[ s_{(2^{k})}s_{(2^{n})}=\sum_{\lambda}a_{\lambda}s_{\lambda}.\]
 By Theorem \ref{lrrule}, the coefficient $a_{\lambda}$ is
equal to the number of Littlewood-Richardson tableaux of shape
$\lambda/(2^n)$ and type $(2^k)$. We claim that $a_{\lambda}=0$ if
the diagram of $\lambda$ contains the square $(n+1,3)$; Otherwise,
we get a contradiction to the assumption $n\geq k$ since the column
strictness of Young tableaux requires that there should be at least
$n+1$ distinct numbers in the tableau. Therefore, for a
Littlewood-Richardson tableau $T$ of shape $\lambda/(2^n)$ and type
$(2^{k})$ we can construct a Littlewood-Richardson tableau $T'$ of
shape $\lambda\cup (2)/(2^{n+1})$ and of the same type by moving all
rows of $T$ to the next row except for the first $n$ rows and
inserting two empty squares at the $(n+1)$-th row. Clearly, the
construction of $T'$ is reversible, as illustrated in Figure
\ref{fig-2}. Hence the first formula is verified. The other
identities can be proved based on similar arguments. \qed

\begin{figure}[h,t]
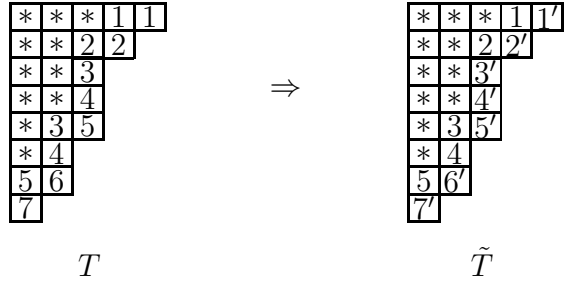

$$
\young{ * & * & * & 1 & 1\cr
        * & * & 2 & 2\cr
        * & * & 3\cr
        * & * & 4\cr
        * & 3 & 5\cr
        * & 4 \cr
        5 & 6 \cr
        7 \cr
        }
\hspace{40pt} \raisebox{47pt}{$\Rightarrow$} \hspace{40pt} \young{
* &
* & * & 1 & 1'\cr
        * & * & 2 & 2'\cr
        * & * & 3'\cr
        * & * & 4'\cr
        * & 3 & 5'\cr
        * & 4 \cr
        5 & 6' \cr
        7' \cr
        }
$$
$$T \hspace{70pt} \hspace{70pt} \tilde{T}$$
\caption{Construct $\tilde{T}$ from $T$}\label{fig-1-1}
\end{figure}

Sometimes it is convenient to regard a tableau $T$ of type
$(2^k,1^l)$ as a semistandard tableau $\tilde{T}$ filled with
distinct numbers in the ordered set
$$\{1<1'<2<2'<\cdots<n<n'<\cdots\}.$$
For this purpose, let $\tilde{T}$ be the tableau such that
$\tilde{T}^{\mathrm{rev}}$ is the word obtained from
$T^{\mathrm{rev}}$ by replacing the first occurrence of $i$ in
$T^{\mathrm{rev}}$ by $i'$ for each $i$ and keeping rest elements
unchanged, as shown in Figure \ref{fig-1-1}.

We also need the following notation to represent a set of partitions
associated with a specified partition. Given a partition $\mu$, let
$$Q_{\mu}(n)=\{\lambda\in \mathrm{Par}(n) : \lambda=\mu \cup (4)^a\cup (3,1)^b\cup (2,2)^c \mbox{ for }a,b,c\in\mathbb{N}\}.$$

\begin{lemm}\label{lemm-odd} Let $m=2k+1$ for some $k\in\mathbb{N}$.
The following statements hold.
\begin{itemize}
\item[(i)]
\begin{equation*}
D_{m,k}^{(1)}=D_{2k+1,k}^{(1)}=s_{(2^{k})}s_{(2^{k})}=\sum_{\lambda\in
Q_{\emptyset}(4k)}s_{\lambda}.
\end{equation*}
\item[(ii)]
\begin{equation*}
D_{m,k+1}^{(1)}=D_{2k+1,k+1}^{(1)}=s_{(2^{k+1})}s_{(2^{k-1})}=\sum_{\lambda\in
Q_{(2,2)}(4k)}s_{\lambda}.
\end{equation*}
\item[(iii)] Let $Q_{1}(n)=Q_{(3,1)}(n)\cup Q_{(2,1,1)}(n)\cup
Q_{(3,3,2)}(n)$. Then
\begin{equation*}
D_{m,k}^{(2)}=D_{2k+1,k}^{(2)}=s_{(2^{k-1},1^2)}s_{(2^{k})}=\sum_{\lambda\in
Q_{1}(4k) }s_{\lambda}.
\end{equation*}

\item[(iv)] Let $Q_2(n)=Q_{(2,1,1)}(n)\cup
Q_{(3,2,2,1)}(n)$. Then
\begin{equation*}
D_{m,k+1}^{(2)}=D_{2k+1,k+1}^{(2)}=s_{(2^{k},1^2)}s_{(2^{k-1})}=\sum_{\lambda\in
Q_2(4k)}s_{\lambda}.
\end{equation*}

\item[(v)] Let $Q_3(n)=Q_{(3,1)}(n)\cup
Q_{(2,2)}(n)\cup Q_{(2,1,1)}(n)\cup Q_{(3,3,2)}(n)$. Then
\begin{equation*}
D_{m,k}^{(3)}=D_{2k+1,k}^{(3)}=s_{(2^{k-1},1)}s_{(2^{k},1)}=\sum_{\lambda\in
Q_3(4k)}a_{\lambda}s_{\lambda},
\end{equation*}
where $a_{\lambda}=2$ if $\lambda\in Q_{(3,2,2,1)}(4k)$, otherwise
$a_{\lambda}=1$.

\item[(vi)] We have
\begin{equation*}
D_{m,k+1}^{(3)}=D_{2k+1,k+1}^{(3)}=s_{(2^{k},1)}s_{(2^{k-1},1)}=\sum_{\lambda\in
Q_3(4k)}a_{\lambda}s_{\lambda},
\end{equation*}
where $a_{\lambda}=2$ if $\lambda\in Q_{(3,2,2,1)}(4k)$, otherwise
$a_{\lambda}=1$.

\end{itemize}
\end{lemm}
\proof
\begin{itemize}
\item[(i)] Use induction on $k$. Clearly, the assertion holds for
$k=0$ since $s_{\emptyset}=1$, and it also holds for $k=1$ by
applying Pieri's rule; see Theorem \ref{pierirule}. From the
Littlewood-Richardson rule it follows that if $s_{\lambda}$ appears
in the Schur expansion of $s_{(2^{k})}s_{(2^{k})}$, then $\lambda$
does not contain any part greater than $4$. So we need to show that
for each Littlewood-Richardson tableau $T$ of shape $\mu/(2^k)$ and
type $(2^k)$, subject to the conditions on the shapes and types,
there are uniquely three Littlewood-Richardson tableaux of type
$(2^{k+1})$, which are $T_1$ of shape $\mu\cup (4)/(2^{k+1})$, $T_2$
of shape $\mu\cup (3,1)/(2^{k+1})$ and $T_3$ of shape $\mu\cup
(2,2)/(2^{k+1})$.

Let $T_1$ be the tableau obtained from $T$ by increasing all numbers
by $1$ and then inserting a four-square row on top of $T$ such that
the rightmost two squares are filled with $1$'s.

Suppose that $T$ has $r$ rows of length greater than $2$, and that
the largest number in the first $r$ rows is $j$ and we set $j=0$ if
$r=0$. Consider the relabeled tableau $\tilde{T}$ corresponding to
$T$. Let $\tilde{T}'$ be the tableau obtained from $\tilde{T}$ by
increasing all numbers below the $r$-th row by $1$ (i.e., changing
$i$ to $i'$ and $i'$ to $i+1$), inserting a three-square row at the
$(r+1)$-th row such that the rightmost square is filled with
$(j+1)'$, and appending a single square row at the bottom filled
with $k+1$. Let $T_2$ be the tableau obtained from $\tilde{T}'$ by
replacing each $i'$ with $i$.

To construct the tableau $T_3$, note that the tableau $T$ does not
contain the square $(k+1,3)$. Consider the numbers  in the first $k$
rows. Let $j_1$ and $j_2$ be the smallest and largest numbers which
appear only once in the first $k$ rows of $T$. Starting with the
tableau $\tilde{T}$, let $\tilde{T}'$ be the tableau obtained from
$\tilde{T}$ by increasing all numbers below the $k$-th row by $2$
(i.e., changing $i$ to $i+1$ and $i'$ to $(i+1)'$), inserting a row
of two empty squares below the $k$-th row, and then inserting a
two-square row filled with $(j_1,(j_2+1)')$ immediately below the
row that has been inserted. If no number appears only once in the
first $k$ rows, consider the largest number $j$ which appears twice
in these rows (taking $j=0$ if no such number exists). Then let
$\tilde{T}'$ be the tableau obtained from $\tilde{T}$ by increasing
all numbers below the $k$-th row by $2$, inserting a row of two
empty squares below the $k$-th row, and then inserting a two-square
row filled with $(j+1,(j+1)')$ immediately below the row just
inserted. Let $T_3$ be the tableau obtained from $\tilde{T}'$ by
replacing each $i'$ with $i$.

Note that if $T$ is a Littlewood-Richardson tableau of shape
$\mu/(2^k)$ and type $(2^k)$, then there exist some nonnegative
integers $r,s,t$ such that the reverse reading word
$\tilde{T}^{\mathrm{rev}}$ is of the form $(w_a,w_b,w_c,w_d)$, where
\begin{eqnarray*}
w_a & = &
1',1,\ldots,r',r\\
w_b & = &
(r+1)',\ldots,(r+s)'\\
w_c & = &
(r+s+1)',(r+1),\ldots,(r+s+t)',(r+s)\\
w_d & = & (r+s+1),\ldots,(r+s+t)
\end{eqnarray*}
and $r+s+t=k$. From $\tilde{T}^{\mathrm{rev}}$ we can write out
${T_1}^{\mathrm{rev}},{T_2}^{\mathrm{rev}},{T_3}^{\mathrm{rev}}$
explicitly according to the above constructions. Now it is easy to
verify that they are lattice permutations.  Figure \ref{fig-3} is an
illustration of the constructions of $T_1,T_2,T_3$.

\begin{figure}[h,t]
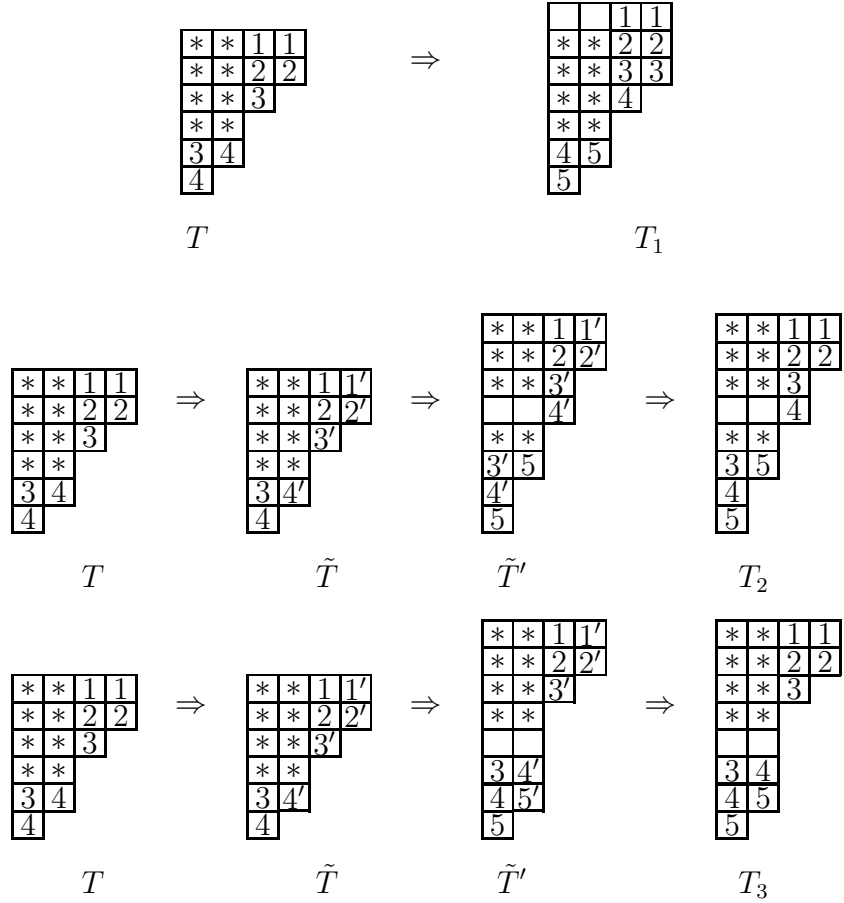

$$
\young{ * & * & 1 & 1\cr
        * & * & 2 & 2\cr
        * & * & 3 \cr
        * & * \cr
        3 & 4 \cr
        4\cr
         }
\hspace{40pt} \raisebox{47pt}{$\Rightarrow$} \hspace{40pt} \young{
  &  & 1 & 1\cr
 * & * & 2 & 2\cr
        * & * & 3 & 3\cr
        * & * & 4 \cr
        * & * \cr
        4 & 5 \cr
        5\cr
         }
$$
$$T \hspace{80pt} \hspace{1pt}\hspace{80pt} T_1$$

$$
\young{ * & * & 1 & 1\cr
        * & * & 2 & 2\cr
        * & * & 3 \cr
        * & * \cr
        3 & 4 \cr
        4\cr
         }
\hspace{15pt} \raisebox{47pt}{$\Rightarrow$} \hspace{15pt} \young{ *
& * & 1 & 1'\cr
        * & * & 2 & 2'\cr
        * & * & 3' \cr
        * & * \cr
        3 & 4' \cr
        4\cr
         }
\hspace{15pt} \raisebox{47pt}{$\Rightarrow$} \hspace{15pt} \young{ *
& * & 1 & 1'\cr
        * & * & 2 & 2'\cr
        * & * & 3' \cr
         &  & 4' \cr
        * & * \cr
        3' & 5 \cr
        4'\cr
        5\cr
         }
\hspace{15pt} \raisebox{47pt}{$\Rightarrow$} \hspace{15pt} \young{ *
& * & 1 & 1\cr
        * & * & 2 & 2\cr
        * & * & 3 \cr
         &  & 4 \cr
        * & * \cr
        3 & 5 \cr
        4\cr
        5\cr
         }
$$
$$T \hspace{80pt} \tilde{T} \hspace{60pt} \tilde{T}' \hspace{80pt} T_2$$
$$
\young{ * & * & 1 & 1\cr
        * & * & 2 & 2\cr
        * & * & 3 \cr
        * & * \cr
        3 & 4 \cr
        4\cr
         }
\hspace{15pt} \raisebox{47pt}{$\Rightarrow$} \hspace{15pt} \young{ *
& * & 1 & 1'\cr
        * & * & 2 & 2'\cr
        * & * & 3' \cr
        * & * \cr
        3 & 4' \cr
        4\cr
         }
\hspace{15pt} \raisebox{47pt}{$\Rightarrow$} \hspace{15pt} \young{ *
& * & 1 & 1'\cr
        * & * & 2 & 2'\cr
        * & * & 3' \cr
        * & * \cr
         &  \cr
        3 & 4' \cr
        4 & 5' \cr
        5\cr
        }
\hspace{15pt} \raisebox{47pt}{$\Rightarrow$} \hspace{15pt} \young{ *
& * & 1 & 1\cr
        * & * & 2 & 2\cr
        * & * & 3 \cr
        * & * \cr
         &  \cr
        3 & 4 \cr
        4 & 5 \cr
        5\cr
        }
$$
$$T \hspace{80pt} \tilde{T} \hspace{60pt} \tilde{T}' \hspace{80pt} T_3$$
\caption{Construction of $T_1,T_2,T_3$ from $T$ of shape
$(4^2,3,2^2,1)/(2^4)$}\label{fig-3}
\end{figure}

On the other hand, it is also necessary to show that for each
Littlewood-Richardson tableau $T'$ of shape $\lambda/(2^{k+1})$ and
type $(2^{k+1})$, we can find a Littlewood-Richardson tableau $T$ of
shape $\mu/(2^k)$ and of type $2^{k}$ such that $\lambda=\mu\cup
(4)$, $\lambda=\mu\cup (3,1)$ or $\lambda=\mu\cup (2,2)$. It is easy
to see that if $\lambda$ contains at least one row of length $4$,
then $T$ can be obtained  from $T'$ by reversing the construction of
$T_1$. If $T'$ has a two-square row fully filled with numbers and
all rows of $T'$ contain at most three squares, then $T$ can be
obtained by reversing the construction of $T_3$. Otherwise, $T'$
contains at least one row of length $1$ and one row of length $3$ in
view of the type of $T'$. In this case, we reverse the construction
of $T_2$ to obtain $T$. Note that $T$ is not uniquely determined by
$T'$. Nevertheless, there exists a unique Littlewood-Richardson
tableau of shape $\lambda/(2^k)$ and type $(2^k)$ if $s_{\lambda}$
appears in the expansion  of $s_{(2^k)}s_{(2^k)}$. Thus,
$s_{(2^k)}s_{(2^k)}$ is multiplicity free.  The proof is completed
by induction.

\item[(ii)] Clearly, the assertion holds for $k=0,1$, and $D_{3,2}^{1}=s_{(2,2)}$.  The
proof is similar to that of (i). Here we consider the
Littlewood-Richardson tableau of shape $\lambda/(2^k)$ and type
$(2^{k-1})$ if $s_{\lambda}$ appears in
$s_{(2^{k+1})}s_{(2^{k-1})}$.

\item[(iii)] Notice that $D_{2k+1,k}^{(2)}=0$ for $k=0$, and
$D_{2k+1,k}^{(2)}=s_{(3,1)}+s_{(2,1,1)}$ for $k=1$. For $k=2$, we
have
\begin{eqnarray*}
D_{2k+1,k}^{(2)}&=&s_{(4,3,1)}+s_{(4,2,1^2)}+s_{(3^2,2)}+s_{(3^2,1^2)}\\
&&+s_{(3,2^2,1)}+s_{(3,2,1^3)}+s_{(2^3,1^2)}.
\end{eqnarray*}
We now use induction on $k$. If $s_{\lambda}$ appears in the
expansion of $D_{2k+1,k}^{(2)}$, then $\lambda$ does not contain the
square $(k+1,3)$, because there exists no Littlewood-Richardson
tableau of shape $\lambda/(2^{k})$ and type $(2^{k-1},1,1)$, or
equivalently, there is no filling of the $(k+1)$-th row satisfying
the lattice permutation condition. Then we can proceed as in the
proof of (i).

\item[(iv)] For $k=0$, it is clear that $D_{2k+1,k+1}^{(2)}=0$. For
$k=1$, we have $D_{2k+1,k+1}^{(2)}=s_{2,1,1}$. For $k=2$, we find
\begin{equation*}
D_{2k+1,k+1}^{(2)}=s_{(4,2,1^2)}+s_{(3,2^2,1)}+s_{(3,2,1^3)}+s_{(2^3,1^2)}.
\end{equation*} Then we use  induction on $k\geq 3$ and
consider Littlewood-Richardson tableaux of shape $\lambda/(2^k,1^2)$
and type $(2^{k-1})$.

\item[(v)] For $k=0$,  $D_{2k+1,k}^{(3)}=0$. For $k=1$, we get
\begin{equation*}
D_{2k+1,k}^{(3)}=s_{(3,1)}+s_{(2^2)}+s_{(2,1^2)}.
\end{equation*}
For $k=2$, we have
\begin{eqnarray*}
D_{2k+1,k}^{(3)}& =&
s_{(4,3,1)}+s_{(4,2^2)}+s_{(4,2,1^2)}+s_{(3^2,2)}\\&&+2s_{(3,2^2,1)}+s_{(3^2,1^2)}+s_{(3,2,1^3)}+s_{(2^4)}+s_{(2^3,1^2)}
\end{eqnarray*}
To use  induction on $k$,  we consider Littlewood-Richardson
tableaux of shape $\lambda/(2^k,1)$ and type $(2^{k-1},1)$. If
$\lambda\in Q_{(3,2,2,1)}(4k)$ there are exactly two such
Littlewood-Richardson tableaux, see Figure \ref{fig-4} for the case
of $\lambda=(4,3^3,2^2,1^3)$. The rest of the proof is similar to
that of (i).

\begin{figure}[h,t]
$$
\young{ * & * & 1 & 1\cr
        * & * & 2 \cr
        * & * & 3 \cr
        * & * & 4 \cr
        *  & * \cr
        *  & 2  \cr
        3\cr
        4\cr
        5\cr
         }
         \hspace{50pt}
\young{ * & * & 1 & 1\cr
        * & * & 2 \cr
        * & * & 3 \cr
        * & * & 4 \cr
        *  & * \cr
        *  & 5 \cr
        2\cr
        3\cr
        4\cr
         }
         $$
\caption{Littlewood-Richardson tableaux of shape
$(4,3^3,2^2,1^3)/(2^5,1)$ and type $(2^4,1)$}\label{fig-4}
\end{figure}

\item[(vi)] It is immediate from (v).
\end{itemize}
This completes the proof of the lemma. \qed

\begin{theo}\label{theo-odd}
Let $m=2k+1$ for some $k\in\mathbb{N}$.
\begin{itemize}
\item[(i)]
We have
\begin{eqnarray}
D_{m,k} & = & s_{(3^k)}s_{(1^k)},\label{eq-oddk}\\
D_{m,k+1}& = &
s_{(4^k)}-s_{(3^k)}s_{(1^k)}-\Delta^{(2)}(s_{(3^k)}s_{1^{(k-2)}}).\label{eq-oddkp}
\end{eqnarray}
\item[(ii)]
For any $0\leq i\leq k-1$, we have
\begin{eqnarray}
D_{m,i} & = & \Delta^{(2)}(D_{m-1,i}),\\
D_{m,m-i} & = & \Delta^{(2)}(D_{m-1,m-1-i}).
\end{eqnarray}
\end{itemize}
\end{theo}

\proof (i) To prove \eqref{eq-oddk}, we need (i), (iii) and (v) of
Lemma \ref{lemm-odd}. If $\lambda\in Q_{(3,2,2,1)}(4k)$, then
$s_{\lambda}$ appears in the expansion of both $D_{m,k}^{(1)}$ and
$D_{m,k}^{(2)}$, and therefore vanishes in $D_{m,k}$. If $\lambda\in
Q_{(3,3,2)}(4k)\cup Q_{(2,1,1)}(4k)$, then $s_{\lambda}$ appears in
both $D_{m,k}^{(2)}$ and $D_{m,k}^{(3)}$, and also vanishes in
$D_{m,k}$. If $\lambda\in Q_{(2,2)}(4k)$ but $\lambda\not\in
Q_{(3,1)}(4k)$, then $s_{\lambda}$ appears in both $D_{m,k}^{(1)}$
and $D_{m,k}^{(3)}$, and also vanishes in $D_{m,k}$. Therefore, for
a term $s_{\lambda}$ which does not vanish in $D_{m,k}$,  the index
partition $\lambda$  belongs to the set $Q_{\emptyset}(4k)$ but $2$
does not appear as a part. By virtue of Pieri's rule, the Schur
functions not vanishing in $D_{m,k}$ coincide with the terms in the
Schur expansion of $s_{(3^k)}s_{(1^k)}$. Similarly, we can prove
\eqref{eq-oddkp} using (ii), (iv) and (vi) of Lemma \ref{lemm-odd}.

(ii) These are direct consequences of Lemma \ref{lemm-first}.

This completes the proof of the theorem.  \qed

Using the same argument in the proof of  Lemma \ref{lemm-odd}, we
can deduce the following expansion formulas when $m$ is even.

\begin{lemm}\label{lemm-even} Let $m=2k$ for $k\in\mathbb{N}$. The
following statements hold.
\begin{itemize}
\item[(i)]
\begin{equation*}
D_{m,k}^{(1)}=D_{2k,k}^{(1)}=s_{(2^{k})}s_{(2^{k-1})}=\sum_{\lambda\in
Q_{(2)}(4k-2)}s_{\lambda}.
\end{equation*}
\item[(ii)]
\begin{equation*}
D_{m,k-1}^{(1)}=D_{2k,k-1}^{(1)}=s_{(2^{k-1})}s_{(2^{k})}=\sum_{\lambda\in
Q_{(2)}(4k-2)}s_{\lambda}.
\end{equation*}
\item[(iii)] Let $R_1(n)=Q_{(1,1)}(n)\cup Q_{(3,3,2,2)}(n)\cup
Q_{(3,2,1)}(n)$. Then
\begin{equation*}
D_{m,k}^{(2)}=D_{2k,k}^{(2)}=s_{(2^{k-1},1^2)}s_{(2^{k-1})}=\sum_{\lambda\in
R_1(4k-2)}s_{\lambda}.
\end{equation*}

\item[(iv)] Let $R_2(n)=Q_{(3,3)}(n)\cup Q_{(3,2,1)}(n)\cup
Q_{(2,2,1,1)}(n)$. Then
\begin{equation*}
D_{m,k-1}^{(2)}=D_{2k,k-1}^{(2)}=s_{(2^{k-2},1^2)}s_{(2^{k})}=\sum_{\lambda\in
R_2(4k-2)}s_{\lambda}.
\end{equation*}

\item[(v)] Let $R_3(n)=Q_{(3,3)}(n)\cup
Q_{(2)}(n)\cup Q_{(1,1)}(n)$. Then
\begin{equation*}
D_{m,k}^{(3)}=D_{2k,k}^{(3)}=s_{(2^{k-1},1)}s_{(2^{k-1},1)}=\sum_{\lambda\in
R_3(4k-2)}a_{\lambda}s_{\lambda},
\end{equation*}
where $a_{\lambda}=2$ if $\lambda\in Q_{(3,2,1)}(4k-2)$, otherwise
$a_{\lambda}=1$.

\item[(vi)] Let
$R_4(n)=Q_{(3,3,2,2)}(n)\cup Q_{(3,2,1)}(n)\cup Q_{(2,2,2)}(n)\cup
Q_{(2,2,1,1)}(n).$
 Then
\begin{equation*}
D_{m,k-1}^{(3)}=s_{(2^{k-2},1)}s_{(2^{k},1)}=\sum_{\lambda\in
R_4(4k-2)}a_{\lambda}s_{\lambda},
\end{equation*}
where $a_{\lambda}=2$ if $\lambda\in Q_{(3,2,2,2,1)}(4k)$, otherwise
$a_{\lambda}=1$.
\end{itemize}
\end{lemm}

In view of Lemmas \ref{lemm-first} and \ref{lemm-even}, we deduce
the following theorem for even $m$. The proof is similar to that of
Theorem \ref{theo-odd} and is omitted.

\begin{theo}\label{theo-even} Let $m=2k$ for some $k\in\mathbb{N}$.
\begin{itemize}
\item[(i)]
We have
\begin{eqnarray}
D_{m,k-1} & = & s_{(3^k)}s_{(1^{k-2})}+\Delta^{(2)}(s_{(3^{k-1})}s_{(1^{k-1})}),\\
D_{m,k} & = & -s_{(3^k)}s_{(1^{k-2})}.\label{eq-evenk}
\end{eqnarray}

\item[(ii)] For any $0\leq i\leq k-2$, we have
\begin{eqnarray}
D_{m,i} & = & \Delta^{(2)}(D_{m-1,i}),\\
D_{m,m-i} & = & \Delta^{(2)}(D_{m-1,m-1-i}),\\
D_{m,m-k+1} & = & \Delta^{(2)}(D_{m-1,m-k}).
\end{eqnarray}
\end{itemize}
\end{theo}

\begin{coro} \label{coro}
Assume $k\geq 1$.
\begin{itemize}
\item[(i)] If $m=2k+1$, then $D_{m,i}$ is $s$-positive for $0\leq i\leq
k$, and $D_{m,i}$ is $s$-negative for $k+1\leq i\leq m-1$.

\item[(ii)]
If $m=2k$, then $D_{m,i}$ is $s$-positive for $0\leq i\leq k-1$, and
$D_{m,i}$ is $s$-negative for $k\leq i\leq m-1$.
\end{itemize}
\end{coro}

\proof Use induction on $m$. It is easy to verify that the result
holds for $k=1$. For $m=2k+1$, we see that $D_{m,k}$ is $s$-positive
and $D_{m,k+1}$ is $s$-negative in view of (i) of Theorem
\ref{theo-odd}. For $0\leq i\leq k-1$, using (ii) of Theorem
\ref{theo-odd} we see that $D_{m,i}=\Delta^{(2)}D_{2k,i}$ is
$s$-positive by induction. Similarly, for $k+2\leq i\leq 2k$ we find
that $D_{m,i}=\Delta^{(2)}D_{2k,i-1}$ is $s$-negative by induction.
For $m=2k$, from (i) of Theorem \ref{theo-even} it follows that
$D_{m,k-1}$ is $s$-positive and $D_{m,k}$ is $s$-negative. For
$0\leq i\leq k-2$, using (ii) of Theorem \ref{theo-even}, by
induction we obtain that $D_{m,i}=\Delta^{(2)}D_{2k-1,i}$ is
$s$-positive. Similarly, for $k+1\leq i\leq 2k-1$, by induction we
deduce that $D_{m,i}=\Delta^{(2)}D_{2k-1,i-1}$ is $s$-negative. \qed

Theorems \ref{theo-odd} and \ref{theo-even} lead to a construction
for the underlying partitions of the Schur expansion. Table
\ref{fig-5} is an illustration.

\begin{table}[p]
\begin{center}
\begin{tabular}{|c|c|}
\hline
& $m=7$\\
\hline
$D_{7,0}$ & $s_{(2^6)}$\\
\hline
$D_{7,1}$ & $s_{(4,2^4)}+s_{(3^2,2^3)}+s_{(3,2^4,1)}$\\
\hline $D_{7,2}$ &
$s_{(3^2,2^2,1^2)}+s_{(4,3^2,2)}+s_{(4^2,2^2)}+s_{(3^3,2,1)}+s_{(4,3,2^2,1)}$\\
\hline $D_{7,3}$ &
$s_{(4,3^2,1^2)}+s_{(3^3,1^3)}+s_{(4^2,3,1)}+s_{(4^3)}$\\
\hline $D_{7,4}$ &
$-s_{(4,3^2,2)}-s_{(4,3^2,1^2)}-s_{(3^3,2,1)}-s_{(3^3,1^3)}-s_{(4^2,3,1)}$\\
\hline $D_{7,5}$ &
$-s_{(3^2,2^3)}-s_{(3^2,2^2,1^2)}-s_{(4,3,2^2,1)}$\\
\hline $D_{7,6}$ &
$-s_{(3,2^4,1)}$\\
\hline $D_{7,7}$ &
$0$\\
\hline
\end{tabular}
\end{center}

\begin{center}
\begin{tabular}{|c|c|}
\hline
& $m=8$\\
\hline
$D_{8,0}$ & $s_{(2^7)}$\\
\hline
$D_{8,1}$ & $s_{(4,2^5)}+s_{(3^2,2^4)}+s_{(3,2^5,1)}$ \\
\hline $D_{8,2}$ &
$s_{(3^2,2^3,1^2)}+s_{(4,3^2,2^2)}+s_{(4^2,2^3)}+s_{(3^3,2^2,1)}+s_{(4,3,2^3,1)}$
\\
\hline $D_{8,3}$ &
$s_{(4,3^2,2,1^2)}+s_{(3^3,2,1^3)}+s_{(4^2,3,2,1)}+s_{(4^3,2)}$
\\
& $+s_{(3^4,1^2)}+s_{(4^2,3^2)}+s_{(4,3^3,1)}$
\\
\hline $D_{8,4}$ & $-s_{(3^4,1^2)}-s_{(4^2,3^2)}-s_{(4,3^3,1)}$
\\
\hline $D_{8,5}$ &
$-s_{(4^2,3,2,1)}-s_{(3^3,2^2,1)}-s_{(3^3,2,1^3)}-s_{(4,3^2,2,1^2)}-s_{(4,3^2,2^2)}$
\\
\hline $D_{8,6}$ &
$-s_{(3^2,2^4)}-s_{(3^2,2^3,1^2)}-s_{(4,3,2^3,1)}$
\\
\hline $D_{8,7}$ & $-s_{(3,2^5,1)}$
\\
\hline $D_{8,8}$ & $0$
\\ \hline
\end{tabular}
\end{center}

\begin{center}
\begin{tabular}{|c|c|}
\hline
& $m=9$\\
\hline
$D_{9,0}$ & $s_{(2^8)}$\\
\hline
$D_{9,1}$ & $s_{(4,2^6)}+s_{(3^2,2^5)}+s_{(3,2^6,1)}$ \\
\hline $D_{9,2}$ &
$s_{(3^2,2^4,1^2)}+s_{(4,3^2,2^3)}+s_{(4^2,2^4)}+s_{(3^3,2^3,1)}+s_{(4,3,2^4,1)}$
\\
\hline $D_{9,3}$ &
$s_{(4,3^2,2^2,1^2)}+s_{(3^3,2^2,1^3)}+s_{(4^2,3,2^2,1)}+s_{(4^3,2^2)}$
\\
 & $+s_{(3^4,2,1^2)}+s_{(4^2,3^2,2)}+s_{(4,3^3,2,1)}$
\\
\hline $D_{9,4}$ &
$s_{(4,3^3,1^3)}+s_{(4^2,3^2,1^2)}+s_{(4^4)}+s_{(4^3,3,1)}+s_{(3^4,1^4)}$
\\
\hline $D_{9,5}$ &
$-s_{(4,3^3,1^3)}-s_{(4^2,3^2,1^2)}-s_{(4^4)}-s_{(4^3,3,1)}-s_{(3^4,1^4)}$
\\
& $-s_{(3^4,2,1^2)}-s_{(4^2,3^2,2)}-s_{(4,3^3,2,1)}$
\\
\hline $D_{9,6}$ &
$-s_{(4^2,3,2^2,1)}-s_{(3^3,2^3,1)}-s_{(3^3,2^2,1^3)}-s_{(4,3^2,2^2,1^2)}-s_{(4,3^2,2^3)}$
\\
\hline $D_{9,7}$ &
$-s_{(3^2,2^5)}-s_{(3^2,2^4,1^2)}-s_{(4,3,2^4,1)}$
\\
\hline $D_{9,8}$ & $-s_{(3,2^6,1)}$
\\
\hline $D_{9,9}$ & $0$
\\ \hline
\end{tabular}
\end{center}
\caption{Schur function expansions of $D_{m,k}$ for
$m=7,8,9$}\label{fig-5}
\end{table}

Given a set $S$ of positive integers, let $\mathrm{Par}_{S}(n)$
denote the set of partitions of $n$ whose parts belong to $S$. We
are now ready to present the following  theorem on Schur positivity.

\begin{theo} \label{theo-s} For any $m\geq 0$, we have
\begin{eqnarray}
&{\sum_{i=0}^m
\left(s_{(2^{i-1})}s_{(2^{m-i})}+s_{(2^{i-2},1^2)}s_{(2^{m-i})}-s_{(2^{i-1},1)}s_{(2^{m-i-1},1)}\right)}&\nonumber\\
&=\sum_{\lambda\in \mathrm{Par}_{\{2,4\}}(2m-2)} s_{\lambda}.&
\label{eq-schursum}
\end{eqnarray}
Consequently, the summation on the left-hand side of the above
identity is $s$-positive.
\end{theo}

Before proving the above theorem, let us give some examples. Taking
$m=3,4,5$ and using the Maple package, we observe that
\begin{eqnarray*}
&\sum_{k=0}^3
\left(s_{(2^{k-1})}s_{(2^{3-k})}+s_{(2^{k-2},1^2)}s_{(2^{3-k})}-s_{(2^{k-1},1)}s_{(2^{3-k-1},1)}\right)&\\
&=s_{(4)}+s_{(2,2)}.&\\ &\sum_{k=0}^4
\left(s_{(2^{k-1})}s_{(2^{4-k})}+s_{(2^{k-2},1^2)}s_{(2^{4-k})}-s_{(2^{k-1},1)}s_{(2^{4-k-1},1)}\right)&\\
&=s_{(4,2)}+s_{(2,2,2)}.& \\
&\sum_{k=0}^5
\left(s_{(2^{k-1})}s_{(2^{5-k})}+s_{(2^{k-2},1^2)}s_{(2^{5-k})}-s_{(2^{k-1},1)}s_{(2^{5-k-1},1)}\right)&\\
&=s_{(4,4)}+s_{(4,2,2)}+s_{(2,2,2,2)}.&
\end{eqnarray*}

\noindent {\it Proof of Theorem \ref{theo-s}.}   By convention, for
$i=0$ or $i=m+1$, it is natural to set
\begin{equation*}
s_{(2^{i-1})}s_{(2^{m-i})}+s_{(2^{i-2},1^2)}s_{(2^{m-i})}=0.
\end{equation*}
Therefore,
\begin{equation*}
{\sum_{i=0}^m
\left(s_{(2^{i-1})}s_{(2^{m-i})}+s_{(2^{i-2},1^2)}s_{(2^{m-i})}-s_{(2^{i-1},1)}s_{(2^{m-i-1},1)}\right)}=\sum_{i=0}^m
D_{m,i}.
\end{equation*}
It suffices to prove that
\begin{equation}\label{eq-mainpf}
\sum_{i=0}^{m+1}
D_{m+1,i}=\left\{\begin{array}{ll}\Delta^{(2)}\left(\sum_{i=0}^m
D_{m,i}\right), & \mbox{if $m=2k-1$}\\[8pt]
s_{(4^k)}+\Delta^{(2)}\left(\sum_{i=0}^m D_{m,i}\right), & \mbox{if
$m=2k$}
\end{array}
\right.
\end{equation}
for $m\geq 0$. The case for $m=0$ is obvious. We now assume $m\geq
1$.

If $m=2k-1$ for some $k\geq 1$, then
$$
\begin{array}{rcl}
\sum_{i=0}^{m+1} D_{m+1,i} & = & \sum_{i=0}^{2k} D_{2k,i} \\[10pt]
                           & = & \sum_{i=0}^{k-2} D_{2k,i} +
                           D_{2k,k-1} + D_{2k,k}+ D_{2k,k+1}+ \sum_{i=0}^{k-2}
                           D_{2k,2k-i}\\[8pt]
                           & = & \sum_{i=0}^{k-2} \Delta^{(2)}(D_{2k-1,i})
                           +
                           \left(s_{(3^k)}s_{(1^{k-2})}+\Delta^{(2)}(s_{(3^{k-1})}s_{(1^{k-1})})\right)\\[8pt]
                           &&+\left(-s_{(3^k)}s_{(1^{k-2})}\right)+\Delta^{(2)}(D_{2k-1,k})\\[8pt]
                           &&+\sum_{i=0}^{k-2}
                           \Delta^{(2)}(D_{2k-1,2k-1-i})\hfill{(\mbox{by Theorem
                           \ref{theo-even}})}\\[10pt]
                           & = & \sum_{i=0}^{k-2} \Delta^{(2)}(D_{2k-1,i})
                           +\Delta^{(2)}(D_{2k-1,k-1})+\Delta^{(2)}(D_{2k-1,k})\\[8pt]
                           &&+\sum_{i=0}^{k-2}
                           \Delta^{(2)}(D_{2k-1,2k-1-i}) \hfill{(\mbox{by
                           \eqref{eq-oddk}})}\\[10pt]
                           &=& \sum_{i=0}^{2k-1} \Delta^{(2)}(D_{2k-1,i})=\Delta^{(2)}\left(\sum_{i=0}^m
                              D_{m,i}\right).
\end{array}
$$
If $m=2k$ for some $k\geq 1$, then
$$
\begin{array}{rcl}
\sum_{i=0}^{m+1} D_{m+1,i} & = & \sum_{i=0}^{2k+1} D_{2k+1,i} \\[10pt]
                           & = & \sum_{i=0}^{k-1} D_{2k+1,i} + D_{2k+1,k}+ D_{2k+1,k+1}+ \sum_{i=0}^{k-1}
                           D_{2k+1,2k+1-i}\\[10pt]
                           & = & \sum_{i=0}^{k-1}
                           \Delta^{(2)}(D_{2k,i})+
                           s_{(3^k)}s_{(1^k)}\\[8pt]
                           &&+\left(s_{(4^k)}-s_{(3^k)}s_{(1^k)}-\Delta^{(2)}(s_{(3^k)}s_{1^{(k-2)}})\right)\\[8pt]
                           &&+ \sum_{i=0}^{k-1}
                           \Delta^{(2)}(D_{2k,2k-i}) \hfill{(\mbox{by Theorem \ref{theo-odd}})}\\[10pt]
                           &=& s_{(4^k)}+\sum_{i=0}^{2k} \Delta^{(2)}(D_{2k,i})\hfill{(\mbox{by \eqref{eq-evenk}})}\\[10pt]
                           &=& s_{(4^k)}+\Delta^{(2)}\left(\sum_{i=0}^m
                               D_{m,i}\right).
\end{array}
$$
 Based on \eqref{eq-mainpf}, we obtain the desired
assertion by induction on $m$. \qed

Now we consider other products of Schur functions, which are
necessary to prove the strong $q$-log-convexity of the Narayana
polynomials.

Given $a,b,r\in\mathbb{N}$ and $0\leq k\leq r$, let
\begin{eqnarray*}
D_1(a,b,k,r) &=&
s_{(2^{k-b-1},1^{b+2-a})}s_{(2^{r-k-1})},\\
D_2(a,b,k,r)&=&
s_{(2^{k-b},1^{b-a})}s_{(2^{r-k-1})},\\
D_3(a,b,k,r)&=&s_{(2^{k-b-1},1^{b+1-a})}s_{(2^{r-k-1},1)}.
\end{eqnarray*}
and let
\begin{equation}
D(a,b,k,r)=D_1(a,b,k,r)+D_2(a,b,k,r)-D_3(a,b,k,r),
\end{equation}
where $s_{(2^{i},1^{j})}=0$ for $i<0$ or $j<0$. It is easy to see
that $D(a,b,r,r)\equiv 0$. For $i=1,2,3$, it is also clear that
$$D_i(a,b,k,r)=D_i(a-1,b-1,k-1,r-1),$$
hence
\begin{equation}
D(a,b,k,r)=D(a-1,b-1,k-1,r-1).
\end{equation}

Some values of $D(a, b, k, r)$ are given in Table \ref{tab-1}.

\begin{table}[p]
\begin{center}
\begin{tabular}{|c|c|}
\hline
& $a=0, b=1, r=8$\\
\hline
$D({a,b,0,r})$ & $0$\\
\hline
$D({a,b,1,r})$ & $s_{(3,2^5)}+s_{(2^6,1)}$\\
\hline $D({a,b,2,r})$ &
$s_{(3^3,2^2)}+s_{(3^2,2^3,1)}+s_{(3,2^4,1^2)}+s_{(4,3,2^3)}+s_{(4,2^4,1)}$\\
\hline $D({a,b,3,r})$ &
$s_{(4,3^2,2,1)}+s_{(4,3,2^2,1^2)}+s_{(3^3,2,1^2)}+s_{(3^2,2^2,1^3)}$\\
& $+s_{(4^2,3,2)}+s_{(4^2,2^2,1)}+s_{(4,3^3)}+s_{(3^4,1)}$\\
\hline $D({a,b,4,r})$ &
$s_{(4^2,3,1^2)}-s_{(4,3^3)}+s_{(4,3^2,1^3)}-s_{(3^4,1)}+s_{(3^3,1^4)}+s_{4^3,1}$\\
\hline $D({a,b,5,r})$ &
$-s_{(4,3^2,2,1)}-s_{(3^3,2^2)}-s_{(3^3,2,1^2)}-s_{(4^2,3,1^2)}-s_{(4,3^2,1^3)}-s_{(3^3,1^4)}$\\
\hline $D({a,b,6,r})$ &
$-s_{(4,3,2^2,1^2)}-s_{(3^2,2^3,1)}-s_{(3^2,2^2,1^3)}$\\
\hline $D({a,b,7,r})$ &
$-s_{(3,2^4,1^2)}$\\
\hline
$D({a,b,8,r})$ & $0$\\
\hline
\end{tabular}
\end{center}

\begin{center}
\begin{tabular}{|c|c|}
\hline
& $a=0, b=1, r=9$\\
\hline
$D({a,b,0,r})$ & $0$\\
\hline
$D({a,b,1,r})$ & $s_{(3,2^6)}+s_{(2^7,1)}$\\
\hline $D({a,b,2,r})$ &
$s_{(3^3,2^3)}+s_{(3^2,2^4,1)}+s_{(3,2^5,1^2)}+s_{(4,3,2^4)}+s_{(4,2^5,1)}$\\
\hline $D({a,b,3,r})$ &
$s_{(4,3^3,2)}+s_{(3^4,2,1)}+s_{4^2,2^3,1)}+s_{(4,3,2^3,1^2)}$\\
& $+s_{(3^2,2^3,1^3)}+s_{(4^2,3,2^2)}+s_{(4,3^2,2^2,1)}+s_{(3^3,2^2,1^2)}$\\
\hline $D({a,b,4,r})$ &
$s_{(4^3,3)}+s_{(4^2,3^2,1)}+s_{(4,3^3,1^2)}+s_{(3^4,1^3)}+s_{(4^3,2,1)}$\\
& $+s_{(4^2,3,2,1^2)}+s_{(4,3^2,2,1^3)}+s_{(3^3,2,1^4)}$\\
\hline $D({a,b,5,r})$ &
$-s_{(4^2,3^2,1)}-s_{(4,3^3,1^2)}-s_{(3^4,1^3)}-s_{(4,3^3,2)}-s_{(3^4,2,1)}$\\
\hline $D({a,b,6,r})$ &
$-s_{(4,3^2,2^2,1)}-s_{(3^3,2^3)}-s_{(3^3,2^2,1^2)}-s_{(4^2,3,2,1^2)}$\\
& $-s_{(4,3^2,2,1^3)}-s_{(3^3,2,1^4)}$\\
\hline $D({a,b,7,r})$ &
$-s_{(4,3,2^3,1^2)}-s_{(3^2,2^4,1)}-s_{(3^2,2^3,1^3)}$\\
\hline $D({a,b,8,r})$ &
$-s_{(3,2^5,1^2)}$\\
\hline
\end{tabular}
\end{center}

\begin{center}
\begin{tabular}{|c|c|}
\hline
& $a=0, b=2, r=10$\\
\hline
$D({a,b,1,r})$ & $0$\\
\hline
$D({a,b,2,r})$ & $s_{(3^2,2^5)}+s_{(3,2^6,1)}+s_{(2^7,1^2)}$\\
\hline $D({a,b,3,r})$ &
$s_{(3^4,2^2)}+s_{(4,3^2,2^3)}+s_{(4,2^5,1^2)}$\\
& $+s_{(3^3,2^3,1)}+s_{(3^2,2^4,1^2)}+s_{(3,2^5,1^3)}+s_{(4,3,2^4,1)}$\\
\hline $D({a,b,4,r})$ & $s_{(4,3^4)}+s_{(4^2,3^2,2)}+s_{(4,3^3,2,1)}+s_{(3^4,2,1^2)}+s_{(4^2,2^3,1^2)}$\\
& $+s_{(4,3,2^3,1^3)}+s_{(3^2,2^3,1^4)}+s_{(4^2,3,2^2,1)}+s_{(4,3^2,2^2,1^2)}+s_{(3^3,2^2,1^3)}$\\
\hline $D({a,b,5,r})$ &
$-s_{(3^5,1)}-s_{(4,3^4)}+s_{(4^3,3,1)}+s_{(4^2,3^2,1^2)}+s_{(4,3^3,1^3)}+s_{(3^4,1^4)}$\\
& $+s_{(4^3,2,1^2)}+s_{(4^2,3,2,1^3)}+s_{(4,3^2,2,1^4)}+s_{(3^3,2,1^5)}$\\
\hline $D({a,b,6,r})$ &
$-s_{(4^2,3^2,1^2)}-s_{(4,3^3,1^3)}-s_{(3^4,1^4)}-s_{(4,3^3,2,1)}-s_{(3^4,2,1^2)}-s_{(3^4,2^2)}$\\
\hline $D({a,b,7,r})$ &
$-s_{(4,3^2,2^2,1^2)}-s_{(3^3,2^2,1^3)}-s_{(3^3,2^3,1)}-s_{(3^3,2,1^5)}$\\
& $-s_{(4,3^2,2,1^4)}-s_{(4^2,3,2,1^3)}$\\
\hline $D({a,b,8,r})$ &
$-s_{(4,3,2^3,1^3)}-s_{(3^2,2^4,1^2)}-s_{(3^2,2^3,1^4)}$\\
\hline $D({a,b,9,r})$ &
$-s_{(3,2^5,1^3)}$\\
\hline
\end{tabular}
\end{center}
\caption{Schur function expansion of $D(a,b,k,r)$ for
$r=8,9,10$}\label{tab-1}
\end{table}

Given a pair $(\lambda,\mu)$ of partitions and a pair $(f_1,f_2)$ of
symmetric functions, we define the product
$\tilde{\Delta}^{\lambda,\mu}(f_1,f_2)$ of $f_1$ and $f_2$ as
follows. Suppose that
\begin{eqnarray}
\Delta^{\lambda}(f_1)&=&\sum_{\nu}a_{\nu}s_{\nu},\\[8pt]
\Delta^{\mu}(f_2)&=&\sum_{\nu}b_{\nu}s_{\nu}.
\end{eqnarray}
Define
\begin{equation}
\tilde{\Delta}^{\lambda,\mu}(f_1,f_2)=\sum_{\nu}\max(a_{\nu},b_{\nu})s_{\nu}.
\end{equation}

\begin{lemm} \label{lemm-two} For any $r\geq k\geq b\geq a\geq 0$ and  $i=1,2,3$,  we have  the following
recurrence relations
\begin{equation}
D_i(a,b,k,r)=\tilde{\Delta}^{(1),(3)}(D_i(a,b-1,k-1,r-1),D_i(a,b-1,k-1,r-2)).
\end{equation}
\end{lemm}

\proof  We first prove that
\begin{eqnarray*}
\lefteqn{ s_{(2^{k-b-1},1^{b+2-a})}s_{(2^{r-k-1})}=}&&
\\&& \tilde{\Delta}^{(1),(3)}(s_{(2^{k-b-1},1^{b+1-a})}s_{(2^{r-k-1})},s_{(2^{k-b-1},1^{b+1-a})}s_{(2^{r-k-2})}).
\end{eqnarray*}
We claim that there exists an injective map between the set of
Littlewood-Richardson tableaux of shape $\mu/(2^{r-k-1})$ and type
$(2^{k-b-1},1^{b+1-a})$ and the set of Littlewood-Richardson
tableaux of shape $\mu\cup (1)/(2^{r-k-1})$ and type
$(2^{k-b-1},1^{b+2-a})$. This means that if $s_{\lambda}$ appears in
the Schur expansion of $s_{(2^{k-b-1},1^{b+1-a})}s_{(2^{r-k-1})}$,
then $s_{\lambda\cup (1)}$ appears in
$s_{(2^{k-b-1},1^{b+2-a})}s_{(2^{r-k-1})}$. This injective map can
be constructed as follows. Given a Littlewood-Richardson tableau $T$
of shape $\mu/(2^{r-k-1})$ and type $(2^{k-b-1},1^{b+1-a})$, let
$T'$ be the tableau obtained from $T$ by appending one row composed
of a single square filled with $k+1-a$. Clearly, $T'$ is a
Littlewood-Richardson tableau of $\lambda\cup (1)/(2^{r-k-1})$ and
type $(2^{k-b-1},1^{b+2-a})$. See the first two tableaux in Figure
\ref{fig-6}.

It will be shown that there exists an injective map between the set
of Littlewood-Richardson tableaux of shape $\mu/(2^{r-k-2})$ and
type $(2^{k-b-1},1^{b+1-a})$ and the set of Littlewood-Richardson
tableaux of shape $\mu\cup (3)/(2^{r-k-1})$ and type
$(2^{k-b-1},1^{b+2-a})$. This means that if $s_{\mu}$ appears in the
expansion of $s_{(2^{k-b-1},1^{b+1-a})}s_{(2^{r-k-2})}$, then
$s_{\mu\cup (3)}$ appears in
$s_{(2^{k-b-1},1^{b+2-a})}s_{(2^{r-k-1})}$. Given a
Littlewood-Richardson tableau $T$ of shape $\mu/(2^{r-k-2})$ and
type $(2^{k-b-1},1^{b+1-a})$, we consider the corresponding tableau
$\tilde{T}$. Suppose that $T$ has $m$ rows of length $4$ (taking
$m=0$ if no such row exists). Let $\tilde{T}'$ be the tableau
obtained from $\tilde{T}$ by inserting one row of three squares at
the $(m+1)$-th row in which the rightmost square is filled with
$(m+1)'$, and then increasing all numbers below the $(m+1)$-th row
by $1$ (i.e., changing $i$ to $i'$ and $i'$ to $i+1$). Let $T'$ be
the tableau obtained from $\tilde{T}'$ by replacing $i'$ with $i$
for each $i$. It is routine to verify that $T'$ is a
Littlewood-Richardson tableau of shape $\mu\cup (3)/(2^{r-k-1})$ and
type $(2^{k-b-1},1^{b+2-a})$, as desired.  See the last four
tableaux in Figure \ref{fig-6}.

Moreover, it remains to prove that any Littlewood-Richardson tableau
$T'$ of shape $\lambda/(2^{r-k-1})$ and type $(2^{k-b-1},1^{b+2-a})$
can be constructed from a Littlewood-Richardson tableau $T$, which
is either of shape $\mu/(2^{r-k-2})$ and type
$(2^{k-b-1},1^{b+1-a})$ with $\lambda=\mu\cup (3)$, or of shape
$\mu/(2^{r-k-1})$ and type $(2^{k-b-1},1^{b+1-a})$ with
$\lambda=\mu\cup (1)$. If $3$ is a part of $\lambda$, then we can
reverse the map in the pervious paragraph to obtain $T$. If $3$ does
not appear as a part of $\lambda$, then the lattice permutation
property requires that $\lambda$ should contain a part of size $1$
and the bottom square should be filled with $k+1-a$. In this case,
let $T$ be the tableau obtained from $T'$ by removing the bottom
row.

Thus we complete the proof of the recurrence of $D_1(a,b,k,r)$, and
the rest can be proved in the same manner. \qed

\begin{figure}[h,t]
$$
\young{ * & * & 1 & 1\cr
        * & * & 2 & 2\cr
        * & * & 3\cr
        * & * & 4\cr
        * & * & 5\cr
        * & * & 6\cr
        * & * \cr
        3 & 7 \cr
        4 & 8 \cr
        }
\hspace{40pt} \raisebox{47pt}{$\Rightarrow$} \hspace{40pt} \young{ *
& * & 1 & 1\cr
        * & * & 2 & 2\cr
        * & * & 3\cr
        * & * & 4\cr
        * & * & 5\cr
        * & * & 6\cr
        * & * \cr
        3 & 7 \cr
        4 & 8 \cr
        9\cr
        }
$$
$$T \hspace{60pt}  \hspace{30pt} \hspace{60pt} T'$$
$$
\young{ * & * & 1 & 1\cr
        * & * & 2 & 2\cr
        * & * & 3\cr
        * & * & 4\cr
        * & * & 5\cr
        * & * \cr
        3 & 6 \cr
        4 & 7 \cr
        8\cr
        }
\hspace{15pt} \raisebox{47pt}{$\Rightarrow$} \hspace{15pt} \young{ *
&
* & 1 & 1'\cr
        * & * & 2 & 2'\cr
        * & * & 3'\cr
        * & * & 4'\cr
        * & * & 5'\cr
        * & * \cr
        3 & 6' \cr
        4 & 7' \cr
        8'\cr
        }
\hspace{15pt} \raisebox{47pt}{$\Rightarrow$} \hspace{15pt} \young{ *
&
* & 1 & 1'\cr
        * & * & 2 & 2'\cr
         &  & 3'\cr
        * & * & 4\cr
        * & * & 5\cr
        * & * & 6\cr
        * & * \cr
        3' & 7 \cr
        4' & 8 \cr
        9\cr
        }
\hspace{15pt} \raisebox{47pt}{$\Rightarrow$} \hspace{15pt} \young{ *
&
* & 1 & 1\cr
        * & * & 2 & 2\cr
        * & * & 3\cr
        * & * & 4\cr
        * & * & 5\cr
        * & * & 6\cr
        * & * \cr
        3 & 7 \cr
        4 & 8 \cr
        9\cr
        }
$$
$$T \hspace{80pt} \tilde{T} \hspace{60pt} \tilde{T}' \hspace{80pt} T'$$
\caption{Two ways to construct $T'$}\label{fig-6}
\end{figure}

\begin{theo} \label{mainconj} For any $b\geq a\geq 0$ and $r\geq 0$,
the symmetric function $\sum_{k=0}^r D(a,b,k,r)$ is $s$-positive.
\end{theo}

\proof We use induction on the difference $b-a$. When $a=b$, note
that
$$\sum_{k=0}^r D(a,b,k,r)=\sum_{k=a}^r D(0,0,k-a,r-a)=\sum_{i=0}^{r-a}D_{r-a,i}.$$
According to Theorem \ref{theo-s}, it is $s$-positive. Now suppose
$b-a\geq 1$. The negative terms of $D(a,b,k,r)$  come from either
${\Delta}^{(1)}(D(a,b-1,k-1,r-1))$ or
${\Delta}^{(3)}(D(a,b-1,k-1,r-2))$ by Lemma \ref{lemm-two}. They
always vanish in $\sum_{k=0}^r D(a,b,k,r)$ since both $\sum_{k=0}^r
D(a,b-1,k-1,r-1)$ and $\sum_{k=0}^r D(a,b-1,k-1,r-2)$ are
$s$-positive by induction. This completes the proof.  \qed

\section{The $q$-Log-convexity}

The main objective of this section is to show that the Narayana
polynomials form a strongly $q$-log-convex sequence. This is a
stronger version of the conjecture of Liu and Wang.

\begin{theo} The Narayana polynomials $N_n(q)$ form a strongly $q$-log-convex
sequence.
\end{theo}
\proof Note that, for $0\leq k\leq n$, we have
\begin{equation}
N(n,k)=N_{q}(n,k)|_{q=1}=s_{(2^k)}(1^{n-1})=\ps_{n-1}^1\left(s_{(2^k)}\right).
\end{equation} For $k> n$, we see that
$N(n,k)=0=\ps_{n-1}^1\left(s_{(2^k)}\right)$.

For any $m\geq n\geq 1$ and $r\geq 0$, the coefficient of $q^r$ in
$N_{m+1}(q)N_{n-1}(q)$ equals
\begin{equation}
C_1=\sum_{k=0}^r \ps_{m}^1\left(s_{(2^k)}\right)
\ps_{n-2}^1\left(s_{(2^{r-k})}\right),
\end{equation}
and the coefficient of $q^r$ in $N_{m}(q)N_{n}(q)$ equals
\begin{equation}
C_2=\sum_{k=0}^r \ps_{m-1}^1\left(s_{(2^k)}\right)
\ps_{n-1}^1\left(s_{(2^{r-k})}\right).
\end{equation}

According to Lemma \ref{lemm-general}, we have
\begin{eqnarray*}
\ps_m^1\left(s_{(2^k)}\right) & = & \sum_{0\leq a\leq b\leq m-n+2}
\ps_{n-2}^1(s_{(2^{k-b},1^{b-a})})\ps_{m-n+2}^1(s_{(2^a,1^{b-a})}),\\[8pt]
\ps_{m-1}^1\left(s_{(2^k)}\right) & = & \sum_{0\leq a\leq b\leq
m-n+1}
\ps_{n-2}^1(s_{(2^{k-b},1^{b-a})})\ps_{m-n+1}^1(s_{(2^a,1^{b-a})}),\\[8pt]
\ps_{n-1}^1\left(s_{(2^{r-k})}\right) & = &
\ps_{n-2}^1\left(s_{(2^{r-k})}+s_{(2^{r-k-1},1)}+s_{(2^{r-k-1})}\right).
\end{eqnarray*}

By expansion we obtain that
$$
\begin{array}{rl}
&C_1-C_2=\\[8pt]
&\hskip 10pt \sum_{k=0}^r\sum_{0\leq a\leq b\leq
m-n+2}\ps_{m-n+2}^1(s_{(2^a,1^{b-a})})\ps_{n-2}^1(s_{(2^{k-b},1^{b-a})}s_{(2^{r-k})})
\\[8pt]
&-\sum_{k=0}^r\sum_{0\leq a\leq b\leq
m-n+1}\ps_{m-n+1}^1(s_{(2^a,1^{b-a})})\ps_{n-2}^1(s_{(2^{k-b},1^{b-a})}s_{(2^{r-k})})\\[8pt]
&-\sum_{k=0}^r\sum_{0\leq a\leq b\leq
m-n+1}\ps_{m-n+1}^1(s_{(2^a,1^{b-a})})\ps_{n-2}^1(s_{(2^{k-b},1^{b-a})}s_{(2^{r-k-1},1)})\\[8pt]
&-\sum_{k=0}^r\sum_{0\leq a\leq b\leq
m-n+1}\ps_{m-n+1}^1(s_{(2^a,1^{b-a})})\ps_{n-2}^1(s_{(2^{k-b},1^{b-a})}s_{(2^{r-k-1})
}).
\end{array}
$$

To simplify the notation, let $d=m-n+1$. Note that
$$
\begin{array}{rcl}
\ps_{d+1}^1(s_{(2^a,1^{b-a})})&=&\ps_{d}^1(s_{(2^a,1^{b-a})})+\ps_{d}^1(s_{(2^a,1^{b-a-1})})\\[8pt]
&&+\ps_{d}^1(s_{(2^{a-1},1^{b-a})})+\ps_{d}^1(s_{(2^{a-1},1^{b-a+1})}).
\end{array}
$$
Therefore, the double summation
$$\sum_{k=0}^r\sum_{0\leq a\leq b\leq
d+1}\ps_{d+1}^1(s_{(2^a,1^{b-a})})s_{(2^{k-b},1^{b-a})}s_{(2^{r-k})}$$
can be divided into four parts
$$
\begin{array}{rcl}
A1 & = & \sum_{k=0}^r\sum_{0\leq a\leq b\leq
d+1}\ps_{d}^1(s_{(2^a,1^{b-a})})s_{(2^{k-b},1^{b-a})}s_{(2^{r-k})}\\[8pt]
A2 & = & \sum_{k=0}^r\sum_{0\leq a\leq b\leq
d+1}\ps_{d}^1(s_{(2^a,1^{b-a-1})})s_{(2^{k-b},1^{b-a})}s_{(2^{r-k})}\\[8pt]
A3 & = & \sum_{k=0}^r\sum_{0\leq a\leq b\leq
d+1}\ps_{d}^1(s_{(2^{a-1},1^{b-a})})s_{(2^{k-b},1^{b-a})}s_{(2^{r-k})}\\[8pt]
A4 & = & \sum_{k=0}^r\sum_{0\leq a\leq b\leq
d+1}\ps_{d}^1(s_{(2^{a-1},1^{b-a+1})})s_{(2^{k-b},1^{b-a})}s_{(2^{r-k})}.
\end{array}
$$

Let
$$
\begin{array}{rcl}
B1 & = & \sum_{k=0}^r\sum_{0\leq a\leq b\leq
d}\ps_{d}^1(s_{(2^a,1^{b-a})})s_{(2^{k-b},1^{b-a})}s_{(2^{r-k})},\\[8pt]
B2 & = & \sum_{k=0}^r\sum_{0\leq a\leq b\leq
d}\ps_{d}^1(s_{(2^a,1^{b-a})})s_{(2^{k-b},1^{b-a})}s_{(2^{r-k-1},1)},\\[8pt]
B3 & = & \sum_{k=0}^r\sum_{0\leq a\leq b\leq
d}\ps_{d}^1(s_{(2^a,1^{b-a})})s_{(2^{k-b},1^{b-a})}s_{(2^{r-k-1}) }.
\end{array}
$$

The equality $A1=B1$ holds because
$$A1=B1+\sum_{k=0}^r\sum_{0\leq a\leq
d+1}\ps_{d}^1(s_{(2^a,1^{d+1-a})})s_{(2^{k-d-1},1^{d+1-a})}s_{(2^{r-k})},$$
but $\ps_{d}^1(s_{(2^a,1^{d+1-a})})\equiv 0$.

We also have the equality $A3=B3$ since
$$
\begin{array}{rcl}
A3 & = & \sum_{k=0}^r\sum_{0\leq a\leq b\leq
d+1}\ps_{d}^1(s_{(2^{a-1},1^{b-a})})s_{(2^{k-b},1^{b-a})}s_{(2^{r-k})}\\[8pt]
 & = & \sum_{k=0}^r\sum_{1\leq a\leq b\leq
d+1}\ps_{d}^1(s_{(2^{a-1},1^{b-a})})s_{(2^{k-b},1^{b-a})}s_{(2^{r-k})}\\[8pt]
 & = & \sum_{k=0}^r\sum_{0\leq a\leq b\leq
d}\ps_{d}^1(s_{(2^{a},1^{b-a})})s_{(2^{k-b-1},1^{b-a})}s_{(2^{r-k})}\\[8pt]
 & = & \sum_{k=1}^r\sum_{0\leq a\leq b\leq
d}\ps_{d}^1(s_{(2^{a},1^{b-a})})s_{(2^{k-b-1},1^{b-a})}s_{(2^{r-k})}\\[8pt]
 & = & \sum_{k=0}^{r-1}\sum_{0\leq a\leq b\leq
d}\ps_{d}^1(s_{(2^{a},1^{b-a})})s_{(2^{k-b},1^{b-a})}s_{(2^{r-k-1})}\\[8pt]
& = & \sum_{k=0}^r\sum_{0\leq a\leq b\leq
d}\ps_{d}^1(s_{(2^a,1^{b-a})})s_{(2^{k-b},1^{b-a})}s_{(2^{r-k-1})
}\\[8pt]
& = & B3.
\end{array}
$$

Moreover, we have
$$
\begin{array}{rcl}
A2 & = & \sum_{k=0}^r\sum_{0\leq a\leq b\leq
d+1}\ps_{d}^1(s_{(2^a,1^{b-a-1})})s_{(2^{k-b},1^{b-a})}s_{(2^{r-k})}\\[8pt]
 & = & \sum_{k=0}^r\sum_{0\leq a<b\leq
d+1}\ps_{d}^1(s_{(2^a,1^{b-a-1})})s_{(2^{k-b},1^{b-a})}s_{(2^{r-k})}\\[8pt]
 & = & \sum_{k=0}^r\sum_{0\leq a\leq b\leq
d}\ps_{d}^1(s_{(2^a,1^{b-a})})s_{(2^{k-b-1},1^{b+1-a})}s_{(2^{r-k})}\\[8pt]
 & = & \sum_{k=0}^r\sum_{0\leq a< b\leq
d}\ps_{d}^1(s_{(2^a,1^{b-a})})s_{(2^{k-b-1},1^{b+1-a})}s_{(2^{r-k})}\\[8pt]
 & & +\sum_{k=0}^r\sum_{0\leq a \leq
d}\ps_{d}^1(s_{(2^a)})s_{(2^{k-a-1},1)}s_{(2^{r-k})}\\[8pt]
 & = & \sum_{k=1}^r\sum_{0\leq a< b\leq
d}\ps_{d}^1(s_{(2^a,1^{b-a})})s_{(2^{k-b-1},1^{b+1-a})}s_{(2^{r-k})}\\[8pt]
 & & +\sum_{k=0}^r\sum_{0\leq a \leq
d}\ps_{d}^1(s_{(2^a)})s_{(2^{k-a-1},1)}s_{(2^{r-k})}\\[8pt]
 & = & \sum_{k=1}^r\sum_{0\leq a\leq b\leq
d-1}\ps_{d}^1(s_{(2^a,1^{b+1-a})})s_{(2^{k-b-2},1^{b+2-a})}s_{(2^{r-k})}\\[8pt]
 & & +\sum_{k=0}^r\sum_{0\leq a \leq
d}\ps_{d}^1(s_{(2^a)})s_{(2^{k-a-1},1)}s_{(2^{r-k})}\\[8pt]
 & = & \sum_{k=0}^{r-1}\sum_{0\leq a\leq b\leq
d-1}\ps_{d}^1(s_{(2^a,1^{b+1-a})})s_{(2^{k-b-1},1^{b+2-a})}s_{(2^{r-k-1})}\\[8pt]
 & & +\sum_{k=0}^r\sum_{0\leq a \leq
d}\ps_{d}^1(s_{(2^a)})s_{(2^{k-a-1},1)}s_{(2^{r-k})}
\end{array}
$$
and
$$
\begin{array}{rcl}
A4 & = & \sum_{k=0}^r\sum_{0\leq a\leq b\leq
d+1}\ps_{d}^1(s_{(2^{a-1},1^{b-a+1})})s_{(2^{k-b},1^{b-a})}s_{(2^{r-k})}\\[8pt]
 & = & \sum_{k=0}^r\sum_{1\leq a\leq b\leq
d}\ps_{d}^1(s_{(2^{a-1},1^{b-a+1})})s_{(2^{k-b},1^{b-a})}s_{(2^{r-k})}\\[8pt]
 & = & \sum_{k=1}^r\sum_{1\leq a\leq b\leq
d}\ps_{d}^1(s_{(2^{a-1},1^{b-a+1})})s_{(2^{k-b},1^{b-a})}s_{(2^{r-k})}\\[8pt]
 & = & \sum_{k=1}^r\sum_{0\leq a\leq b\leq
d-1}\ps_{d}^1(s_{(2^{a},1^{b+1-a})})s_{(2^{k-b-1},1^{b-a})}s_{(2^{r-k})}\\[8pt]
 & = & \sum_{k=0}^{r-1}\sum_{0\leq a\leq b\leq
d-1}\ps_{d}^1(s_{(2^{a},1^{b+1-a})})s_{(2^{k-b},1^{b-a})}s_{(2^{r-k-1})}
\end{array}
$$
and
$$
\begin{array}{rcl}
B2 & = & \sum_{k=0}^r\sum_{0\leq a\leq b\leq
d}\ps_{d}^1(s_{(2^a,1^{b-a})})s_{(2^{k-b},1^{b-a})}s_{(2^{r-k-1},1)}\\[8pt]
 & = & \sum_{k=0}^r\sum_{0\leq a< b\leq
d}\ps_{d}^1(s_{(2^a,1^{b-a})})s_{(2^{k-b},1^{b-a})}s_{(2^{r-k-1},1)}\\[8pt]
&&+\sum_{k=0}^r\sum_{0\leq a\leq
d}\ps_{d}^1(s_{(2^a)})s_{(2^{k-a})}s_{(2^{r-k-1},1)}\\[8pt]
 & = & \sum_{k=0}^r\sum_{0\leq a< b\leq
d}\ps_{d}^1(s_{(2^a,1^{b-a})})s_{(2^{k-b},1^{b-a})}s_{(2^{r-k-1},1)}\\[8pt]
 & & +\sum_{k=0}^r\sum_{0\leq a \leq
d}\ps_{d}^1(s_{(2^a)})s_{(2^{k-a-1},1)}s_{(2^{r-k})}\\[8pt]
 & = & \sum_{k=0}^{r-1}\sum_{0\leq a\leq b\leq
d-1}\ps_{d}^1(s_{(2^a,1^{b+1-a})})s_{(2^{k-b-1},1^{b+1-a})}s_{(2^{r-k-1},1)}\\[8pt]
 & & +\sum_{k=0}^r\sum_{0\leq a \leq
d}\ps_{d}^1(s_{(2^a)})s_{(2^{k-a-1},1)}s_{(2^{r-k})}.
\end{array}
$$
Therefore,
$$
\begin{array}{rcl}
C_1-C_2 & = & \ps_{n-2}^1((A_1+A_2+A_3+A_4)-(B_1+B_2+B_3))\\[8pt]
        & = & \ps_{n-2}^1(A_2+A_4-B_2)\\[8pt]
        & = & \ps_{n-2}^1\left(\sum_{0\leq a\leq b\leq
d-1}\ps_{d}^{1}(s_{(2^a,1^{b+1-a})})\sum_{k=0}^{r} D(a,b,k,r)\right)
\end{array}
$$
From Theorem \ref{mainconj} we deduce that
$$\sum_{0\leq a\leq b\leq
d-1}\ps_{d}^{1}(s_{(2^a,1^{b+1-a})})\sum_{k=0}^{r} D(a,b,k,r)$$ is
$s$-positive, hence $C_1-C_2$ is nonnegative, as desired. \qed

As a corollary, we are led to an affirmative answer to Conjecture
\ref{mainprob}.

\begin{coro}\label{coro-nara} The Narayana polynomials $N_n(q)$ form a $q$-log-convex
sequence.
\end{coro}

\textbf{Remark.} Butler and Flanigan \cite{butfla2007} defined a
different  $q$-analogue of log-convexity. In their definition, a
sequence of polynomials $(f_k(q))_{k\geq 0}$ is called
$q$-log-convex if
$$f_{m-1}(q)f_{n+1}(q)-q^{n-m+1}f_{m}(q)f_{n}(q)$$
has nonnegative coefficients for $n\geq m\geq 1$. They proved that
the $q$-Catalan numbers of Carlitz and Riordan \cite{cr20} form a
$q$-log-convex sequence. However, the Narayana polynomial sequence
$(N_{n}(q))_{n\geq 0}$ is not $q$-log-convex  by the definition of
Butler and Flanigan.

\section{The Narayana transformation}

In \cite{liuwan2006} Liu and Wang studied several log-convexity
preserving transformations, and they also realized the connection
between the $q$-log-convexity and the linear transformations
preserving the log-convexity. They conjectured that if the sequence
$(a_k)_{k\geq 0}$ of positive real numbers is log-convex then the
sequence $$b_n=\sum_{k=0}^n N(n,k) a_k, \quad {n\geq 0}$$ is also
log-convex. In this section we will provide a proof of this
conjecture. We first give two  lemmas.

For any $n\geq 1$ and $0\leq r\leq 2n$, we define the following
polynomials in  $x$ with integer coefficients:
\begin{eqnarray*}
f_1(x) & = & (n+1)(n-x+1)(n-x)^2(n-x-1),\\[8pt]
f_2(x) & = & (n+1)(n-(r-x)+1)(n-(r-x))^2(n-(r-x)-1),\\[8pt]
f_3(x)  & = & (n-1)(n-x)(n-x+1)(n-(r-x))(n-(r-x)+1).
\end{eqnarray*}
Let $$f(x)=f_1(x)+f_2(x)-2f_3(x).$$
\begin{lemm}\label{lem1} For fixed integers $n\geq 1$ and $0\leq r< 2n$, the polynomial $f(x)$ is
monotone decreasing in $x$ on the interval $(-\infty,\frac{r}{2}]$.
\end{lemm}
\proof Taking the derivative $f'(x)$ of $f(x)$ with respect to $x$,
we obtain that
 $$f'(x)=2(2x-r)g(x),$$
 where
 $$g(x)=4x^2-4xr-2n+r+2r^2-5nr-8n^2r-2+2nr^2+6n^2+8n^3.$$
Note that the discriminant of the quadratic polynomial $g(x)$ equals
\begin{eqnarray}\label{eq-q1}
&(-4r)^2-16(-2n+r+2r^2-5nr-8n^2r-2+2nr^2+6n^2+8n^3)\nonumber&\\[8pt]
&=16(-r^2+2n-r+5nr+8n^2r+2-2nr^2-6n^2-8n^3).&
\end{eqnarray}
Let us consider the following polynomial
$$g_1(y)=-y^2+2n-y+5ny+8n^2y+2-2ny^2-6n^2-8n^3$$
in $y$ on the interval $(-\infty,2n)$. The derivative of $g_1(y)$
with respect to $y$ is
$$g_1'(y)=-2y-1+5n+8n^2-4ny=(4n+2)(2n-y)+n-1.$$
Therefore, $g_1'(y)> 0$ for $y\in (-\infty,2n)$. Then for any $0\leq
r<2n$ and $n\geq 1$ we have
$$g_1(r)\leq g_1(2n-1)=-3n+2<0.$$
This implies that $g(x)>0$ and $f'(x)=2(2x-r)g(x)<0$ for
$x\in(-\infty,\frac{r}{2})$. Therefore, $f(x)$ is monotone
decreasing on the interval $(-\infty,\frac{r}{2}]$. \qed

\begin{lemm}\label{lemm-con}
For any $n\geq 1$, $0\leq r\leq 2n$ and $0\leq k\leq \lfloor
\frac{r}{2}\rfloor$, let
\begin{eqnarray*}
\alpha(n,r,k)&=&N(n+1,k)N(n-1,r-k)+N(n+1,r-k)N(n-1,k)\\[8pt]
&&-2N(n,r-k)N(n,k).
\end{eqnarray*}
Then, for given $n$ and $r$,  there always exists an integer
$k'=k'(n,r)$ such that $\alpha(n,r,k)\geq 0$ for $k\leq k'$ and
$\alpha(n,r,k)\leq 0$ for $k> k'$.
\end{lemm}

\proof Assume that $n$ and $r$ are given. Clearly, if $k\leq r-n-1$,
then $n\leq (r-k)-1$ and $\alpha(n,r,k)=0$. We only need to
determine the sign of $\alpha(n,r,k)$ for $r-n-1<k\leq \lfloor
\frac{r}{2}\rfloor$.

Note that $N(m,k)=\ps_{m-1}^{1}(s_{(2^k)})$ for any $m\in
\mathbb{N}$. By Lemma \ref{lemm-hook} we find that
$$N(m,k)=\frac{((n-1)(n-2)\cdots(n-r+k))\cdot (n(n-1)\cdots(n-r+k+1))}{k!(k+1)!}.$$
Let
\begin{eqnarray*}
C & = & (n-1)(n-2)^2(n-3)^2\cdots(n-k+2)^2 (n-k+1),\\[8pt]
C' & = & (n-1)(n-2)^2(n-3)^2\cdots(n-(r-k)+2)^2 (n-(r-k)+1).
\end{eqnarray*}
Then we have
\begin{equation*}
\alpha(n,r,k)= \frac{C}{k!(k+1)!}\cdot\frac{C'}{(r-k)!(r-k+1)!}\cdot
f(k).
\end{equation*}

Now let us consider the value of $f(k)$ for fixed $r$. We have the
following three cases.
\begin{itemize}
\item[(i)] When $r=2m+1$ for some $0\leq m<n$, by Lemma \ref{lem1} we have
$$f(0)\geq f(1)\geq \cdots \geq f(m),$$
where
$$f(0)=2(2m+1)(n+1)((4m+1)(n-m)(n-m-1)+m(m+1))\geq 0.$$

\item[(ii)] When $r=2m$ for some $0\leq m <n$, by Lemma \ref{lem1} we have
$$f(0)\geq f(1)\geq \cdots \geq f(m),$$
where
$$f(0)=4m(n+1)((4m-1)(n-m)^2+m(m-1)\geq 0.$$

\item[(iii)] When $r=2n$, we have
$$f(k)=4(n-k+1)(n-k-1)(n-k)^2.$$
Therefore, $f(k)>0$ for any $k< n$ and $f(n)=0$.
\end{itemize}

Notice that there always exists an integer $k'$ such that $f(k)\geq
0$ for $k\leq k'$ and $f(k)\leq 0$ for $k> k'$. Because both $C$ and
$C'$ are nonnegative, we reach the desired conclusion. \qed

\begin{theo}\label{trans}
If the sequence $(a_k)_{k\geq 0}$ of positive real numbers is
log-convex, then the sequence
\[ b_n=\sum_{k=0}^n N(n,k) a_k, \quad n\geq
0 \]
 is
log-convex.
\end{theo}

In general, the Narayana transformation does not preserve the
log-convexity, and the condition that $(a_k)_{k\geq 0}$ is a
positive sequence is necessary for the above theorem. For example,
if we take $a_k=(-1)^k$ for $k\geq 0$, then it is easy to see that
$(a_k)_{k\geq 0}$ is log-convex, but $(b_n)_{n\geq 0}$ is not
log-convex.

\noindent \textit{Proof of Theorem \ref{trans}.} For any $n,r,k\geq
0$, let
$$\alpha'(n,r,k)=\left\{\begin{array}{ll}
\alpha(n,r,k)/2, & \mbox{if $r$ is even and $k=r/2$},\\[8pt]
\alpha(n,r,k), & \mbox{otherwise}.
\end{array}\right.$$
Note that for $n\geq 1$
$$b_{n-1}b_{n+1}-b_n^2=\sum_{r=0}^{2n}\left(\sum_{k=0}^{\lfloor \frac{r}{2}\rfloor}\alpha'(n,r,k)a_ka_{r-k}\right)$$
and
$$N_{n-1}(q)N_{n+1}(q)-N_n(q)^2=\sum_{r=0}^{2n}\left(\sum_{k=0}^{\lfloor \frac{r}{2}\rfloor}\alpha'(n,r,k)\right)q^r.$$
By Corollary \ref{coro-nara}, we see that
$$\sum_{k=0}^{\lfloor
\frac{r}{2}\rfloor}\alpha'(n,r,k)\geq 0$$  for any $r\geq 0$. Since
the sequence $(a_k)_{k\geq 0}$ is a log-convex sequence of positive
real numbers, we obtain that
$$a_0a_r\geq a_1a_{r-1}\geq a_2a_{r-2}\geq \cdots.$$
 Lemma \ref{lemm-con} implies that there exists an integer $k'=k'(n,r)$ such
that
$$
\sum_{k=0}^{\lfloor \frac{r}{2}\rfloor}\alpha'(n,r,k)a_ka_{r-k}\geq
\sum_{k=0}^{\lfloor
\frac{r}{2}\rfloor}\alpha'(n,r,k)a_{k'}a_{r-k'}\geq 0.
$$
Therefore, $(b_n)_{n\geq 0}$ is log-convex. \qed

\section{The $q$-log-concavity}

This section is devoted to the $q$-log-concavity of the $q$-Narayana
numbers $N_q(n,k)$ for given $n$ or $k$.  First we apply
Br\"and$\mathrm{\acute{e}}$n's formula \eqref{formula-b} to express
the $q$-Narayana numbers in terms of specializations of Schur
functions. This formulation enables us to reduce the
$q$-log-concavity of the $q$-Narayana numbers to the  Schur
positivity of some  differences between the products of Schur
functions indexed by two-column shapes. Notice that much work has
been done on the Schur positivity of the differences of products of
Schur functions, see, for example, Bergeron, Biagioli and Rosas
\cite{bebiro2005}, Fomin, Fulton, Li and Poon \cite{fflp2005} and
Okounkov \cite{okounk1997}.

We now proceed to prove the $q$-log-concavity of $q$-Narayana
numbers $N_q(n,k)$ for fixed $n$.

\begin{theo}\label{logcave1}
Given a positive integer $n$, the sequence $(N_q(n,k))_{k\geq 0}$ is
strongly $q$-log-concave.
\end{theo}

\proof Using \eqref{formula-b}, for any $k\geq l\geq 1$, we get
\begin{equation*}
N_q(n,k)N_q(n,l)-N_q(n,k+1)N_q(n,l-1)
=s_{(2^{k})}s_{(2^{l})}-s_{(2^{k+1})}s_{(2^{l-1})},
\end{equation*}
where each Schur function on the righthand side is over the variable
set $\{q,q^2,\ldots,q^{n-1}\}$. Using induction on $k-l$, we can
show that the symmetric function
$s_{(2^{k})}s_{(2^{l})}-s_{(2^{k+1})}s_{(2^{l-1})}$ is $s$-positive.
Clearly, this statement is true for $k=l$, because by (i) and (ii)
of Lemma \ref{lemm-odd}, we have
\begin{eqnarray*}
s_{(2^{k})}s_{(2^{k})}- s_{(2^{k+1})}s_{(2^{k-1})}  =
\sum_{\lambda\in Q_{\emptyset}(4k)}s_{\lambda}- \sum_{\lambda\in
Q_{(2,2)}(4k)}s_{\lambda}=\sum_{a=0}^ks_{(4^a,3^{k-a},1^{k-a})}.
\end{eqnarray*}

For $k>l$, by \eqref{eq-s-1} we have
\begin{eqnarray*}
s_{(2^{k})}s_{(2^{l})} & = & \Delta^{(2)}(s_{(2^{k-1})}s_{(2^{l})}),\\[5pt]
s_{(2^{k+1})}s_{(2^{l-1})}  & = & \Delta^{(2)}
(s_{(2^{k})}s_{(2^{l-1})}).
\end{eqnarray*}
It follows that
\begin{eqnarray*}
s_{(2^{k})}s_{(2^{l})} - s_{(2^{k+1})}s_{(2^{l-1})}  =
\Delta^{(2)}(s_{(2^{k-1})}s_{(2^{l})}-s_{(2^{k})}s_{(2^{l-1})}).
\end{eqnarray*}
By induction, $s_{(2^{k-1})}s_{(2^{l})}-s_{(2^{k})}s_{(2^{l-1})}$ is
$s$-positive, so is $s_{(2^{k})}s_{(2^{l})} -
s_{(2^{k+1})}s_{(2^{l-1})}$.
The Schur positivity of the above difference was also shown by
Bergeron and McNamara \cite[Remark 7.2]{bermcn2004}, and Kleber
\cite{kleber2001} gave a  proof for the case of $k=l$. In view of
the variable set for symmetric functions, we see that the difference
$N_q(n,k)N_q(n,l)-N_q(n,k+1)N_q(n,l-1)$ has nonnegative coefficients
as a polynomial of $q$. This completes the proof.\qed

Next we will consider the $q$-log-concavity of the $q$-Narayana
numbers $N_q(n,k)$ for fixed $k$. We will use a result due to Lam,
Postnikov and Pylyavaskyy \cite{lapopy2007}. Given two partitions
$\lambda=(\lambda_1,\lambda_2,\ldots)$ and
$\mu=(\mu_1,\mu_2,\ldots)$, let
\begin{eqnarray*}
\lambda \vee \mu & = &
(\max(\lambda_1,\mu_1),\max(\lambda_2,\mu_2),\ldots),\\
\lambda \wedge \mu & = &
(\min(\lambda_1,\mu_1),\min(\lambda_2,\mu_2),\ldots).
\end{eqnarray*}
For two skew partitions $\lambda/\mu$ and $\nu/\rho$, we define
\begin{eqnarray*}
(\lambda/\mu) \vee (\nu/\rho) & = & (\lambda\vee \nu)/(\mu \vee
\rho),\\[3pt]
(\lambda/\mu) \wedge (\nu/\rho) & = & (\lambda\wedge \nu)/(\mu
\wedge \rho).
\end{eqnarray*}

The following assertion was conjectured by Lam and Pylyavaskyy
\cite{lampyl2007} and proved by Lam, Postnikov and Pylyavaskyy
\cite{lapopy2007}. We will be interested in two special cases of
this fact.

\begin{theo}[{\cite[Theorem 5]{lapopy2007}}]\label{sect-theo}
For any two skew partitions $\lambda/\mu$ and $\nu/\rho$, the
difference
$$s_{(\lambda/\mu) \vee (\nu/\rho)}s_{(\lambda/\mu) \wedge
(\nu/\rho)}-s_{\lambda/\mu}s_{\nu/\rho}$$ is $s$-positive.
\end{theo}

In particular, we will need the following special cases.

\begin{coro}\label{sect-6-coro}
Let $k$ be a positive integer. If $I,J$ are partitions with
$I\subseteq (2^{k-1})$ and $J\subseteq (2^{k-1},1)$, then both
\begin{equation}\label{eq-1}
s_{(2^{k-1})}s_{(2^{k})/I}- s_{(2^{k-1})/I}s_{(2^{k})}
\end{equation}
and
\begin{equation}\label{eq-2}
s_{(2^{k-1},1)}s_{(2^{k})/J}- s_{(2^{k-1},1)/J}s_{(2^{k})}
\end{equation}
are $s$-positive.
\end{coro}
\proof For \eqref{eq-1}, take $\lambda=(2^{k-1}), \mu=I, \nu=(2^k)$
and $\rho=\emptyset$ in Theorem \ref{sect-theo}. For \eqref{eq-2},
take $\lambda=(2^{k-1},1), \mu=J, \nu=(2^k)$ and $\rho=\emptyset$.
\qed

For any $r\geq 1$, let
$$X_r=\{q,q^2,\ldots,q^{r-1}\}, \quad
X_r^{-1}=\{q^{-1},q^{-2},\ldots,q^{-(r-1)}\}.$$ The following
relations are crucial for the proof of the $q$-log-concavity of the
$q$-Narayana numbers $N_q(n,k)$ for given $k$.

\begin{lemm}\label{sect-6-lemm} For any $m\geq n\geq 1$ and $k\geq 1$, we have
\begin{eqnarray}
\lefteqn{q^{n-1}s_{(2^{k-1},1)}(X_{n-1})s_{(2^{k})}(X_m)-
q^{m}s_{(2^{k-1},1)}(X_{m})s_{(2^{k})}(X_{n-1})}  &  &\nonumber \\[5pt]
& = & q^{k-1}\left(s_{(2^{k-1},1)}(X_{n-1})s_{(2^{k})}(X_m)-
s_{(2^{k-1},1)}(X_{m})s_{(2^{k})}(X_{n-1})\right)\label{eq-inv-1}
\end{eqnarray}
and \begin{eqnarray}
\lefteqn{q^{2(n-1)}s_{(2^{k-1})}(X_{n-1})s_{(2^{k})}(X_m)-
q^{2m}s_{(2^{k-1})}(X_{m})s_{(2^{k})}(X_{n-1})} \nonumber \\[5pt]
 & = & q^{2k(m+n-1)}\left(s_{(2^{k-1})}(X_{n-1}^{-1})s_{(2^{k})}(X_m^{-1})-
s_{(2^{k-1})}(X_{m}^{-1})s_{(2^{k})}(X_{n-1}^{-1})\right). \hskip
1cm\label{eq-inv-2}
\end{eqnarray}
\end{lemm}
\proof We will adopt the following notation for $q$-series in the
proof. For  indeterminates $a,a_1,\cdots,a_s$ and integer $r\geq 0$,
let
$$
\begin{array}{ccc}
(a;q)_r & = & (1-a)(1-aq)\cdots (1-aq^{r-1}),\\[5pt]
(a_1,a_2,\cdots,a_s;q)_r & = & (a_1;q)_r(a_2;q)_r\cdots (a_s;q)_r.
\end{array}
$$
 By Lemma \ref{lemm-hook}, we
have
\begin{eqnarray*}
s_{(2^{k-1},1)}(X_{n-1}) & = &
s_{(2^{k-1},1)}(q,q^2,\cdots,q^{n-2})\\[5pt]
 & = &
\frac{q^{k^2}(q^{n-k-1};q)_k(q^{n-k+1};q)_{k-1}}{(1-q)(q;q)_{k-1}
(q^3;q)_{k-1}}
\end{eqnarray*}
and
\begin{eqnarray*}
s_{(2^{k})}(X_{n}) & = &
s_{(2^{k})}(q,q^2,\cdots,q^{n-1})\\[5pt]
 & = &
 \frac{q^{k(k+1)}(q^{n-k};q)_{k}(q^{n-k+1};q)_{k}}{(q;q)_{k}(q^2;q)_{k}}.
\end{eqnarray*}
Therefore, the left hand side of \eqref{eq-inv-1} equals
\begin{eqnarray*}
\lefteqn{\frac{q^{2k^2+k+n-1}(q^{n-k+1};q)_{k-1}(q^{n-k-1},q^{m-k},q^{m-k+1};q)_{k}}{(1-q)(q,q^3;q)_{k-1}
(q,q^2;q)_{k}}}&&\\[5pt]
&-&
\frac{q^{2k^2+k+m}(q^{m-k+2};q)_{k-1}(q^{m-k},q^{n-k-1},q^{n-k};q)_{k}}{(1-q)(q,q^3;q)_{k-1}
(q,q^2;q)_{k}}\\[5pt]
&=&\frac{q^{2k^2+k+n-1}(1-q^{m-n+1})(q^{m-k+2},q^{n-k+1};q)_{k-1}(q^{m-k},q^{n-k-1};q)_{k}}
{(1-q)(q,q^3;q)_{k-1} (q,q^2;q)_{k}}
\end{eqnarray*}
and the difference $s_{(2^{k-1},1)}(X_{n-1})s_{(2^{k})}(X_m)-
s_{(2^{k-1},1)}(X_{m})s_{(2^{k})}(X_{n-1})$ equals
\begin{eqnarray*}
\lefteqn{\frac{q^{2k^2+k}(q^{n-k+1};q)_{k-1}(q^{n-k-1},q^{m-k},q^{m-k+1};q)_{k}}{(1-q)(q,q^3;q)_{k-1}
(q,q^2;q)_{k}}}&&\\[5pt]
&-&
\frac{q^{2k^2+k}(q^{m-k+2};q)_{k-1}(q^{m-k},q^{n-k-1},q^{n-k};q)_{k}}{(1-q)(q,q^3;q)_{k-1}
(q,q^2;q)_{k}}\\[5pt]
&=&\frac{q^{2k^2+n}(1-q^{m-n+1})(q^{m-k+2},q^{n-k+1};q)_{k-1}(q^{m-k},q^{n-k-1};q)_{k}}
{(1-q)(q,q^3;q)_{k-1} (q,q^2;q)_{k}}.
\end{eqnarray*}
Comparing the above two identities, we arrive at \eqref{eq-inv-1}.

Next we prove the second identity. The left hand side of
\eqref{eq-inv-2} equals
\begin{eqnarray*}
\lefteqn{\frac{q^{2(n+k^2-1)}(q^{n-k},q^{n-k+1};q)_{k-1}(q^{m-k},q^{m-k+1};q)_{k}}{(q,q^2;q)_{k-1}
(q,q^2;q)_{k}}}&&\\[5pt]
&-&
\frac{q^{2(m+k^2)}(q^{m-k+1},q^{m-k+2};q)_{k-1}(q^{n-k-1},q^{n-k};q)_{k}}{(q,q^2;q)_{k-1}
(q,q^2;q)_{k}}\\[5pt]
&=&\frac{f(q)(q^{n-k},q^{n-k+1},q^{m-k+1},q^{m-k+2};q)_{k-1}}{(q,q^2;q)_{k-1}
(q,q^2;q)_{k}},
\end{eqnarray*}
where
$$f(q)=q^{2k^2-k-2}(q^{m+1}-q^n)(q^{m+n+1}+q^{m+n}-q^{m+k+1}-q^{n+k}).$$

The difference $s_{(2^{k-1})}(X_{n-1})s_{(2^{k})}(X_m)-
s_{(2^{k-1})}(X_{m})s_{(2^{k})}(X_{n-1})$ equals
\begin{eqnarray*}
\lefteqn{\frac{q^{2k^2}(q^{n-k},q^{n-k+1};q)_{k-1}(q^{m-k},q^{m-k+1};q)_{k}}{(q,q^2;q)_{k-1}
(q,q^2;q)_{k}}}&&\\[5pt]
&-&
\frac{q^{2k^2}(q^{m-k+1},q^{m-k+2};q)_{k-1}(q^{n-k-1},q^{n-k};q)_{k}}{(q,q^2;q)_{k-1}
(q,q^2;q)_{k}}\\[5pt]
&=&\frac{g(q)(q^{n-k},q^{n-k+1},q^{m-k+1},q^{m-k+2};q)_{k-1}}{(q,q^2;q)_{k-1}
(q,q^2;q)_{k}},
\end{eqnarray*}
where
$$g(q)=q^{2k^2-2k-1}(q^{m+1}-q^n)(q^{m+1}+q^{n}-q^{k+1}-q^{k}).$$

It is routine to verify that $g(q^{-1})=q^{2k+1-4k^2-2m-2n}f(q)$.
Then \eqref{eq-inv-2} follows from the fact that
$(1-q^{-r})=-q^{-r}(1-q^r)$ for any $r$.
 \qed

Now we are ready to prove the $q$-log-concavity of the $q$-Narayana
numbers $(N_q(n,k))_{n\geq 0}$ for given $k$.

\begin{theo}\label{logcave1}
Given a positive integer $k$, the sequence $(N_q(n,k))_{n\geq 0}$ is
strongly $q$-log-concave.
\end{theo}

\proof For any $m\geq n\geq 1$, let
$$A_{m,n}(q)  =
N_q(m,k)N_q(n,k)-N_q(m+1,k)N_q(n-1,k).$$ By \eqref{formula-b}, we
have
$$
\begin{array}{rcl}
A_{m,n}(q) & = &
s_{(2^{k})}(X_m)s_{(2^{k})}(X_n)-s_{(2^{k})}(X_{m+1})s_{(2^{k})}(X_{n-1}).
\end{array}
$$

Applying \eqref{eq-expansion} to $s_{(2^{k})}(X_n)$ and
$s_{(2^{k})}(X_{m+1})$, the above $A_{m,n}(q)$ equals
\begin{eqnarray*}
\lefteqn{
s_{(2^{k})}(X_m)\left(s_{(2^{k})}(X_{n-1})+q^{n-1}s_{(2^{k-1},1)}(X_{n-1})+q^{2(n-1)}s_{(2^{k-1})}(X_{n-1})\right)}&&\\[8pt]
&&-\left(s_{(2^{k})}(X_{m})+q^{m}s_{(2^{k-1},1)}(X_{m})+q^{2m}s_{(2^{k-1})}(X_{m})\right)s_{(2^{k})}(X_{n-1})\\[8pt]
&=&\left( q^{n-1}s_{(2^{k-1},1)}(X_{n-1})s_{(2^{k})}(X_m)-
q^{m}s_{(2^{k-1},1)}(X_{m})s_{(2^{k})}(X_{n-1})\right)\\[8pt]
&&+ \left( q^{2(n-1)}s_{(2^{k-1})}(X_{n-1})s_{(2^{k})}(X_m)-
q^{2m}s_{(2^{k-1})}(X_{m})s_{(2^{k})}(X_{n-1})\right).
\end{eqnarray*}

By Lemma \ref{sect-6-lemm}, we obtain that $A_{m,n}(q)$ equals
$$
\begin{array}{ll}
&\hskip -6mm q^{k-1}\left(s_{(2^{k-1},1)}(X_{n-1})s_{(2^{k})}(X_m)-
s_{(2^{k-1},1)}(X_{m})s_{(2^{k})}(X_{n-1})\right)\\[8pt]
&+
q^{2k(m+n-1)}\left(s_{(2^{k-1})}(X_{n-1}^{-1})s_{(2^{k})}(X_m^{-1})-
s_{(2^{k-1})}(X_{m}^{-1})s_{(2^{k})}(X_{n-1}^{-1})\right)\\[8pt]
&\hskip -4mm =
q^{k-1}s_{(2^{k-1},1)}(X_{n-1})s_{(2^{k})}(Z)\\[8pt]
&+q^{k-1}\sum_{J\subseteq
(2^{k-1},1)}s_J(Z)\left(s_{(2^{k-1},1)}s_{(2^{k})/J}-
s_{(2^{k-1},1)/J}s_{(2^{k})}\right)(X_{n-1})\\[8pt]
&+q^{2k(m+n-1)}s_{(2^{k-1})}(X_{n-1}^{-1})s_{(2^{k})}(Z^{-1})\\[8pt]
&+q^{2k(m+n-1)}s_{(2^{k-1})}(X_{n-1}^{-1})s_{(2^{k-1},1)}(Z^{-1})s_{(1)}(X_{n-1}^{-1})\\[8pt]
&+q^{2k(m+n-1)}\sum_{I\subseteq
(2^{k-1})}s_I(Z)\left(s_{(2^{k-1})}s_{(2^{k})/I}-
s_{(2^{k-1})/I}s_{(2^{k})}\right)(X_{n-1}^{-1}),
\end{array}
$$
where $Z=\{q^{n-1},\ldots,q^{m-1}\}$ and
$Z^{-1}=\{q^{1-n},\ldots,q^{1-m}\}$. Applying Corollary
\ref{sect-6-coro}, we complete the proof.  \qed

\vskip 8pt

\noindent {\bf Acknowledgments.} This work was supported by  the 973
Project, the PCSIRT Project of the Ministry of Education, the
Ministry of Science and Technology, and the National Science
Foundation of China.

\end{document}